%% file: x3njn1111.tex
\documentclass[12pt]{article}
\usepackage{times,a4wide,amsmath,latexsym,amssymb}
\newtheorem{theorem}{Theorem}[section]
\newtheorem{lemma}[theorem]{Lemma}
\newtheorem{proposition}[theorem]{Proposition}

\newtheorem{rem}[theorem]{Remark}
\newtheorem{cor}[theorem]{Corollary}
\newtheorem{prop}[theorem]{Proposition}

\newcommand{\8}{\infty}

\newcommand{\kla}{\left ( }
\newcommand{\mer}{\right ) }

\newcommand{\for}{\begin{eqnarray*}}
\newcommand{\mel}{\end{eqnarray*}}

\newcommand{\kl}{\pl \le \pl}
\newcommand{\gl}{\pl \ge \pl}

\newcommand{\lel}{\pl = \pl}

\newcommand{\bN}{{\mathbb N}}
\newcommand{\bC}{{\mathbb C}}

\newcommand{\nz}{{\rm  I\! N}}
\newcommand{\nen}{n \in \nz}

\newcommand{\cz}{{\Bbb C}}

\newcommand{\ten}{\otimes}

\newcommand{\p}{\hspace{.05cm}}
\newcommand{\pl}{\hspace{.1cm}}

\newcommand{\hz}{\vspace{0.5cm}}

\newcommand{\qed}{\hspace*{\fill}$\Box$\hz\pagebreak[1]}

\newcommand{\Om}{\Omega}
\newcommand{\om}{\omega}

\newcommand{\al}{\alpha}
\newcommand{\si}{\sigma}

\newcommand{\eps}{\varepsilon}

\newcommand{\cU}{{\cal U}}
\newcommand{\cN}{{\cal N}}

\newcommand{\cL}{{\cal L}}
\newcommand{\cO}{{\cal O}}
\newcommand{\cC}{{\cal C}}
\newcommand{\cS}{{\cal S}}

\renewcommand{\L}{{\cal L}}

\newcommand{\E}{{\cal E}}

\newcommand{\M}{{\cal M}}
\newcommand{\bproof}{\noindent{\bf Proof: }}
\newcommand{\eproof}{\hfill $\Box$\\}

\newcommand{\N}{{\cal N}}

\newcommand{\U}{{\cal U}}
\newcommand{\K}{{\cal K}}

\newcommand{\noo}{\left \|}
\newcommand{\rrm}{\right \|}

\newcommand{\bet}{\left |}
\newcommand{\rag}{\right |}

\newcommand{\intt}{\int\limits}
\newcommand{\summ}{\sum\limits}

\newcommand{\prodd}{\prod\nolimits}

\newcommand{\ez}{{\rm I\!E}}

\renewcommand{\span}{{\rm span}}
\newcommand{\re}{\begin{rem}\rm}
\begin{document}

\title{Rosenthal operator spaces}
\author{M.\ Junge\footnote{Supported by NSF grant DMS--0301116 and
DMS 05-56120}, N.J.\ Nielsen\footnote{Supported by the Danish
Natural Science Research Council, grant 21020436.}, and  T.\
Oikhberg\footnote{Supported by NSF grant DMS--0500957}} \footnotetext{2000 {\em Mathematics Subject
Classification:} 46B20,
  46L07, 46L52.}
\date{}

\maketitle

\input njnpart1new11

\input njnpart2new11

\input rosmatrixnew11
\input uncomp2.tex

\input refnew11.tex

\end{document}

%% file: njnpart1new11.tex
\begin{abstract}

\noindent In 1969 Lindenstrauss and Rosenthal showed that if a Banach
space is isomorphic to a complemented subspace of an $L_p$-space, then it
is either a $\cL_p$-space or isomorphic to a Hilbert space. This is
the motivation of this paper where we study non--Hilbertian complemented
 operator subspaces of non commutative $L_p$-spaces and show that this class
is much richer than in the commutative case. We investigate the local properties
of some new classes of operator spaces for every $2<p< \infty$ which can be considered as operator space analogues of the Rosenthal sequence spaces from Banach
space theory, constructed in 1970. Under the usual conditions on the defining sequence $\sigma$ we prove that most of these spaces are operator
$\cL_p$-spaces, not completely isomorphic to previously known such spaces.
However it turns out that some column and row versions of our spaces are not
operator $\cL_p$-spaces and have a rather complicated local structure
which implies that the Lindenstrauss--Rosenthal alternative does not
carry over to the non-commutative case.

\end{abstract}

\section*{Introduction}

In 1970 Rosenthal \cite{R} constructed new examples of $\cL_p$--spaces 
for every $2\le p <\infty$ using probabilistic methods now famous as the
Rosenthal inequalities. These methods were later used by Bourgain, Rosenthal 
and Schechtman \cite{BRS} to construct an uncountable family of mutually 
non-isomorphic $\cL_p$--spaces.

\indent In the framework of operator spaces a
theory of operator $\cL_p$-spaces, called $\cO\cL_p$-spaces, is now being
developed, see e.g. \cite{ER1} and \cite{JNRX}. These are spaces where the
operator space structure of the finite dimensional subspaces is determined by
a system of finite dimensional non commutative $L_p$-spaces. If in a given
space these $L_p$-spaces can be chosen to be completely complemented, the
space is called a $\cC\cO\cL_p$-space. If they can be chosen to be $S_p^n$'s
($S_p$ denotes the Schatten $p$-class), then the space is called an
$\cO\cS_p$-space and a $\cC\cO\cS_p$-space if the $S_p^n$'s can be chosen completely complemented. In the present paper we consider some operator space
analogues of the Rosenthal sequence spaces, sequence spaces as well as
matricial analogues.

For a given $2<p< \infty$ and a given strictly positive sequence
$\sigma = (\sigma_n)$ we construct three families of operator spaces, a
sequence space family consisting of spaces called $X_p(\sigma)$,
$X_{p,r_p}(\sigma)$ and $X_{p,c_p}(\sigma)$, and two families of matricial operator spaces.
All the spaces are mutually non-completely isomorphic as operator spaces, but the spaces in each family are isomorphic to each other as Banach spaces; the
three sequence spaces are actually Banach space isomorphic to the original Rosenthal
sequence space. One of our main results states that if
$2<p<\infty$, $\sigma_n \to 0$, and
$\sum_{n=1}^{\infty}\sigma_n^{\frac{2p}{2p-2}}=\infty$, then
$X_{p,c_p}(\sigma)$ is completely complemented in a non commutative $L_p$-space
and contains $\ell_p$ cb-complemented. However $X_{p,c_p}(\sigma)$ is not an
$\cO\cL_p$-space. Similarly for $X_{p,r_p}(\sigma)$. This shows that the Lindenstrauss-Rosenthal
alternative \cite{LR} does not carry over to the non commutative case.

\indent We now wish to discuss the arrangement of this paper in greater detail.
In Section 1 we construct our spaces, investigate their basic properties and
prove among other things that under the above conditions on
$\sigma$ the three sequence spaces are unique up to complete isomorphisms
(in analogy with Rosenthal's result). In Section 2 we
make a detailed investigation of the local structure of the spaces
$X_p(\sigma)$, $X_{p,c_p}(\sigma)$ and $X_{p,r_p}(\sigma)$ and prove that
$X_p(\sigma)$ is an $\cO\cL_p$-space while $X_{p,r_p}(\sigma)$ and
$X_{p,c_p}(\sigma)$ are not. We also show that some combinations of the
different spaces cannot be paved with local pieces of each other. This
implies that a general structure theory for completely complemented
non-Hilbertian subspaces of non commutative $L_p$-spaces is out of reach for the moment (see e.g. Proposition \ref{propultra} and Remark \ref{remultra}).
Section 3 is devoted to the study of the
matricial spaces and we show that they are all $\cO\cS_p$-spaces and
prove that the space $Y_p(\sigma)$ is cb-complemented in $L_p(\mathcal
R)$ ($\mathcal R$
the hyperfinite type $II_1$ factor) while the space $Z_p(\sigma)$ does
not cb-embed into $L_p(\mathcal R)$. In section 4 we prove that certain 
$\cO\cL_p$--spaces contain cb-uncomplemented copies of themselves.

\setcounter{section}{-1}
\section{Notation and preliminaries}
\label{sec0}
\setcounter{equation}{0}

In this paper we shall use the notation and terminology commonly used in the
theory of operator algebras, operator spaces and Banach space theory as it appears in \cite{ER2}, \cite{JNRX}, \cite{LT1}, \cite{LT2}, \cite{P1} and \cite{Ta}.

If $H$ is a Hilbert space, we let $B(H)$ denote the space of all bounded operators on $H$ and for every $n \in \bN$ we let $M_n$ denote the space of all
$n \times n$-matrices of complex numbers, i.e. $M_n = B(\ell_{2}^n)$. If $X$ is
a subspace of some $B(H)$ and $n \in \bN$, then $M_n(X)$ denotes the space of
all $n \times n$ matrices with $X$-valued entries which we in the natural manner consider as a
subspace of $B(\ell_2 ^n (X))$. An {\em operator space} $X$ is a norm closed
subspace
of some $B(H)$ equipped with the distinguised matrix norm inherited by the
spaces $M_n(X)$, $n \in \bN$. An abstract matrix norm characterization of operator spaces was given by Ruan, see e.g. \cite{ER2}.

If $X$ and $Y$ are operator spaces, then a linear operator $T: X \rightarrow Y$ is
called {\em completely bounded} (in short cb-bounded) if the corresponding
linear maps $T_n : M_n(X) \rightarrow M_n(Y)$ are uniformly bounded in $n$,
i.e.
\begin{equation*}
\|T\|_{cb} = sup\|T_n\|< \infty
\end{equation*}

The space of all completely bounded operators from $X$ to $Y$ will be denoted by
$CB(X,Y)$.

It follows from \cite{ER2} that a linear functional on an operator space $X$ is
bounded if and only if it is cb-bounded and the cb-norm of it coincides with the operator norm of it. This defines an operator structure on $X^*$ so that
isometrically we have $M_n(X^*)= CB(X,M_n)$ for all $n \in \bN$.

An operator is a {\em complete contraction}, respectively a {\em complete
isometry}, or a {\em complete quotient} if $\|T\|_{cb} \le 1$, respectively
if each $T_n$ is an isometry, or a quotient map. An operator $T$ is called
a {\em complete isomorphism} (in short a {\em cb-isomorphism}) if it is a
completely bounded linear isomorphism with a completely bounded linear
inverse. If $X$ and $Y$ are cb-isomorphic operator spaces we put
\begin{equation*}
d_{cb}(X,Y) = inf\{\|T\|_{cb} \|T^{-1}\|_{cb} \mid \mbox{T is a cb-isomorphism
from X to Y} \}
\end{equation*}
which is called the {\em completely bounded Banach--Mazur distance} (in short
the {\em cb-distance}) between $X$ and $Y$.

In the sequel we let $S_{\infty} \subseteq B(\ell_2)$ denote the subspace of all
compact operators on $\ell_2$ (hence an operator space in a natural manner).
If $1\le p < \infty$, then the {\em Schatten class} $S_p$ is defined to be the
space of all compact operators $T$ on $\ell_2$ for which $tr(|T|)^p < \infty$
equipped with the norm
\begin{equation}
\label{eq01}
\|T\|_{S_p} = (tr(|T|^p))^{\frac{1}{p}} \quad \mbox{for all $T \in S_p$}
\end{equation}

If $n \in \bN$ and $p$ is as above, $S_p ^n$ denotes the space of all operators
on $\ell_2 ^n$ equipped with the norm defined in (\ref{eq01}). If also
$m\in \bN$, then $S_p ^{n,m}$ denotes the subspace of $S_p$ consisting of
those elements which correspond to matrices $(a_{ij})$ where $a_{ij} = 0$
unless $i\le n$ and $j\le m$.

From trace duality it easily follows that $S_{\infty}^* = S_1$ and hence as
a dual space $S_1$ has a natural operator structure as defined above. It is
wellknown that $S_p$ can be obtained by by complex interpolation
\begin{equation*}
S_p = [S_{\infty}, S_1]_{\frac{1}{p}}
\end{equation*}

Pisier proved in \cite{P1} that
\begin{equation*}
M_n(S_p) = [M_n(S_{\infty}), M_n(S_1)]_{\frac{1}{p}}
\end{equation*}

defines matrix norms on $S_p$ which satisfy Ruan's matrix norm characterization
of operator spaces and this is called the {\em natural operator space
structure} of $S_p$ which we shall always use in the sequel.

Let $e_{ij}$ denote the element of $B(\ell_2)$ corresponding to the matrix
with coefficients equal to one at the $i,j$ entry and zero elsewhere. If
$1\le p \le \infty$, we define the operator subspaces $C_p$ and $R_p$ of
$S_p$ by
\begin{eqnarray*}
C_p = \overline{span}\{e_{i1} \mid i\in \bN \}\\
R_p = \overline{span} \{e_{1j} \mid j\in \bN \}.
\end{eqnarray*}

As Banach spaces these spaces are both isometric to $\ell_2$, but it follows from Pisier \cite{P1} that they are not cb-isomorphic as operator spaces.

If $1 \le p \le \infty$, then we put ${\K}_p = (\sum_n^{\infty}
S_p^n)_p$; ${\K}_p$ is clearly an operator space in a canonical manner.

If $H$ is an operator Hilbert space, i.e. an operator space which as a Banach space is isometric to a Hilbert space, then we put $H^c = CB(\bC,H)$ and
$H^r=CB(H,\bC)$ and if $1<p<\infty$, then we let $H^{c_p}= [H^c,H^r]_{\frac1p}$
and $H^{r_p}=[H^r,H^c]_{\frac1p}$.

If $E$ is an operator space and $1 \le p \le \infty$, it is possible to
define $S_p[E]$ ($S_p$ with values in$E$) as the completion of $S_p\otimes
E$ under a certain operator space norm; we refer to \cite[chapter
1]{P1} for the details. In particular we shall often use the following
proposition proved by Pisier \cite[Lemma 1.7, see also Propositions 2.3, 2.4 and Remark
2.5]{P1}.

\begin{proposition}
\label{prop01}

Let $E$ and $F$ be operator spaces. A linear map $T: E \to F$ is cb-bounded 
if and only $\sup_{n\in \bN} \|Id_{S_p ^n}\otimes T: S_p^n[E]\to S_p^n[F]\|< \infty$. 
In the affirmative we have 
$\|T\|_{cb} = \sup_{n\in \bN} \|Id_{S_p ^n}\otimes T \|$.
\end{proposition}

The norms in $S_p[R_p]$ and $S_p[C_p]$ were computed by Pisier in
\cite[page 108]{P1} and since we are going to use this frequently in the sequel
we state it in a proposition.

\begin{proposition}
\label{prop02}

If $(x_k)_{k=1}^n \subseteq S_p$, then
\begin{equation}
\label{pr}
\|\sum_{k=1}^n x_k \otimes e_{1k}\|_{S_p[R_p]} =
\|(\sum_{k=1}^n x_k x_k^*)^{\frac{1}{2}}\|_{S_p}
\end{equation}
and
\begin{equation}
\label{pc}
\|\sum_{k=1}^n x_k \otimes e_{k1}\|_{S_p[C_p]} =
\|(\sum_{k=1}^n x_k^* x_k)^{\frac{1}{2}}\|_{S_p}
\end{equation}
\end{proposition}

If $X$ is a subspace of $S_p$ and $E$ is an operator space, then we let
$X[E]$ denote the closure of $E\otimes X$ in $S_p[E]$.

Let $A$ be a von Neumann algebra with a normal semifinite faithful trace $\tau$
(i.e. $A$ is semifinite). The ideal
\begin{equation*}
m(\tau)= \{\sum_{k=1}^n x_k y_k \mid n\in \bN, \quad \sum_{k=1} ^n[\tau(y_k ^* y_k)+ \tau(x_k^* x_k)]< \infty \}
\end{equation*}
is called the definition ideal of $\tau$ on which there is a unique linear
extension $\tau: m(\tau) \rightarrow \bC$ so that $\tau(xy)= \tau(yx)$ for all
$x,y \in m(\tau)$ (see e.g \cite{Ta}). If $1\le p < \infty$, then we put
\begin{equation*}
\|x\| = \tau((x^* x)^{\frac{p}{2}})^{\frac{1}{p}} \quad \mbox{for all $x\in
m(\tau)$}
\end{equation*}
which is readily seen to be a norm on $m(\tau)$. We define $L_p(A,\tau)$ to be
the completion of $m(\tau)$ under this norm. Conventionally we put
$L_{\infty}(A,\tau)= A$. It follows easily that $L_1(A, \tau)^* = A^{op}$
where $A^{op}$ denotes $A$ equipped with the reversed (or opposite) multiplication and hence
$L_1(A,\tau)$ has a natural operator space structure. It can be shown that
the complex interpolation method yields that
\begin{equation*}
L_p(A,\tau) = [A,L_1(A,\tau)]_{\frac{1}{p}}.
\end{equation*}

Pisier \cite{P1} proved that
\begin{equation*}
M_n(L_p(A,\tau))= [M_n(A),M_n(L_1(A,\tau))]_{\frac{1}{p}}
\end{equation*}

defines a natural operator space structure on $L_p(A,\tau)$ which we shall
use in the sequel. If ${\tau}_1$ is another normal semifinite faithful trace
on $A$, then it can easily be shown that $L_p(A,\tau)$ is cb-isometric to
$L_p(A,{\tau}_1)$ and therefore we shall often write $L_p(A)$ instead of
$L_p(A,\tau)$.

If $B$ is von Neumann subalgebra of A so that the restriction of $\tau$ to
$B$ is semifinite again, then it follows from \cite[Proposition 2.36]{Ta} that
there exists a faithful normal projection $E_B$ of $A$ onto $B$ such that
$\tau=\tau \circ E_B$. $E_B$ is called {\em the conditional expectation} of $A$
onto $B$.

An operator space $X$ is called an {\em operator ${\cL}_p$-space} (in short
$\cO{\cL}_p$- space, $1\le p \le \infty$, if there exist a $\lambda \ge 1$ and a
cofinal family $(F_j)_{j\in I}$ of finite dimensional subspaces so that $\bigcup_{j\in I} F_j$ is dense in $X$ and so that for every index $j$ there exists a finite dimensional $C^*$-algebra $A_j$ with
\begin{equation}
\label{eq02}
d_{cb}(L_p(A_j),F_j) \le \lambda.
\end{equation}

In this case we shall also say that $X$ is a $\cO{\cL}_{p,\lambda}$-space. $X$
is called an $\cO{\cS}_{p,\lambda}$-space if we can replace the $L_p(A_j)$'s
in (\ref{eq02}) by $S_p ^{n_j}$'s. $X$ is called a {\em completely
complemented $\cO{\cL}_{p,\lambda}$-space} (in short $\cC\cO{\cL}_{p,\lambda}$-space), if in addition the $F_j$'s can be chosen to be cb-complemented in $X$ by
projections with cb-norms less than or equal to $\lambda$. $\cC \cO {\cS}_{p,\lambda}$-spaces are defined similarly.

If the $L_p(A_j)$'s in (\ref{eq02}) are of the form
$(\oplus_{i=1}^k S_p ^{n(i),m(i)})_p$, then $X$ is called a {\em rectangular
$\cO \cL_p$-space}.

Let $1 \le p \le \infty$. An operator space $X$ is said to have the
{\em $\gamma_p$-approximation property} (in short $\gamma_p$-AP) if
there exists a $\lambda > 0$ and nets $(U_i)$ and $(V_i)$ of finite rank
operators, $U_i\colon X \to S_p$, $V_i \colon S_p \to X$, so that
$\|U_i\|_{cb} \|V_i\|_{cb} \le \lambda$ and $(V_i U_i)$ converges
pointwise to the identity of $X$.

Finally, if $(x_n)$ is a finite or infinite sequence in a Banach space
$X$, we let $[x_n]$ denote the closed linear span of the sequence
$(x_n)$. If $A$ is a set, $|A|$ denotes the cardinality of $A$ and if $X$ and $Y$ are Banach spaces, $X \oplus_p Y$ denotes the direct sum of $X$ and $Y$ equipped with the norm $(\| \cdot \|_X^p + \|\cdot\|_Y^p)^{\frac{1}{p}}$.

\section{The Rosenthal operator spaces and their basic properties}
\label{sec1}
\setcounter{equation}{0}

In this section we shall investigate some operator spaces
which in nature correspond to the $\cL_p$-spaces in Banach space theory
constructed by Rosenthal in \cite{R}.

In the sequel we let $2<p<\infty$, $\frac1p +\frac1{p'} =1$,
$\frac12 = \frac1p +\frac1r$ (i.e. $r= \frac{2p}{p-2}$)
and let $\sigma = (\sigma_n)$ be a sequence of real
numbers with $\sigma_n > 0$ for all $n\in
\bN$.  We denote the unit vector basis of
$\ell_2$ by $(\xi_n)$ and let
$D_{\bf \sigma}$ be the diagonal operator on $\ell_2$ defined by
$D_{\sigma}\xi_n=\sigma_n \xi_n$ for all $n\in \bN$.

Our first space $\tilde{X}_{p}(\sigma)$ is defined to be the space of all
sequences $a= (a_n)$ which satisfies
\begin{equation}
\label{eq1.3}
\sum_{n=1}^{\infty} |a_n|^p < \infty \quad \mbox{and} \quad
\sum_{n=1}^{\infty} |a_n|^2 \sigma_n^2 < \infty.
\end{equation}
equipped with the norm

\begin{equation}
\|a\|= (\sum_{n=1}^{\infty}|a_n|^p +
(\sum_{n=1}^{\infty}|a_n|^2 \sigma_n^2)^{\frac{p}{2}})^{\frac{1}{p}}
\end{equation}

$\tilde{X}_p(\sigma)$ is the classical Rosenthal sequence space (except that he used an equivalent norm) and we
can clearly identify it
with the closed linear span in $S_p\oplus_p S_2$ of the sequence
$\{(e_{nn},\sigma_n e_{nn}) \mid n\in \bN\}$. As an operator space we can
however represent $\tilde{X}_p(\sigma)$ in three different ways. We define the
space $X_{p,c_p}(\sigma)$ to be
the closed linear span of the sequence $\{(e_{nn},\sigma_n e_{n1}) \mid
n\in\bN \}$ in $S_p \oplus C_p$. Similarly we let $X_{p,r_p}(\sigma)$
denote the closed linear span of the sequence $\{(e_{nn},\sigma_n e_{1n}) \mid
n\in \bN \}$ in $S_p \oplus R_p$ and finally we let $X_p(\sigma)$ denote
the closed linear span of the sequence $\{(e_{nn},\sigma_n
e_{n1},\sigma_n e_{1n})\}$ in $S_p \oplus C_p\oplus R_p$.

Since $S_p \oplus C_p \oplus R_p$ is cb-isomorphic to $S_p$ each of the three
above spaces is cb-isomorphic to a subspace of $S_p$. In the sequel we
shall often let $X_{p*}(\sigma)$ denote any of these spaces.

Since we shall often use Proposition \ref{prop01} to check cb-boundedness in this paper it is worthwhile to mention how the norms in $S_p[X_{p,r_p}(\sigma)]$,
$S_p[X_{p,c_p}(\sigma)]$ and $S_p[X_p(\sigma)]$ can be computed. It follows immediately from
Proposition \ref{prop02} that if $(x_k)_{k=1}^n \subseteq S_p$, then
\begin{equation}
\label{row}
\|\sum_{k=1}^n x_k \otimes (e_{kk} \oplus \sigma_k e_{1k}) \|_{S_p[X_{p,r_p}(\sigma)]} =
(\sum_{k=1}^n \|x_k\|^p + \|(\sum_{k=1}^n \sigma_k^2
x_k x_k^*)^{\frac{1}{2}}\|_{S_p}^p)^{\frac{1}{p}}
\end{equation}

\begin{equation}
\label{raekke}
\|\sum_{k=1}^n x_k \otimes (e_{kk} \oplus \sigma_k e_{k1})\|_{S_p[X_{p,c_p}(\sigma)]} =
(\sum_{k=1}^n \|x_k\|^p + \|(\sum_{k=1}^n \sigma_k^2
x_k^* x_k)^{\frac{1}{2}} \|_{S_p}^p)^{\frac{1}{p}}
\end{equation}

and similarly for $S_p[X_p(\sigma)]$.

It follows easily from these formulas and Proposition \ref{prop01} that
though isometric as Banach spaces
these three spaces are not mutually cb-isomorphic as operator spaces.

Throughout the paper we shall often impose at least one of the following
two conditions on $\sigma$:
\begin{equation}
\label{eq1.1}
\liminf_{n\rightarrow \infty}\sigma_n = 0
\end{equation}

\begin{equation}
\label{eq1.2}
\sum_{\sigma_n \le \varepsilon} \sigma_n^r =
\infty \quad \quad \mbox{for all $\varepsilon > 0$}
\end{equation}

It is immediate that if $\sigma_n \rightarrow 0$ and
$\sigma \notin \ell_r$, then (\ref{eq1.1}) and (\ref{eq1.2}) are
satisfied. \eqref{eq1.2}
ensures that the operator $x \rightarrow xD_{\sigma}$ does not act as a
bounded operator from $S_p$ to $S_2$.

It follows from \cite{R} that $\tilde{X}_p(\sigma)$ is an $\cL_p$-space
if and only if (\ref{eq1.1}) is satisfied, and if both (\ref{eq1.1}) and
(\ref{eq1.2}) holds, then $\tilde{X}_p(\sigma)$ is the classical
Rosenthal $\cL_p$-space which is unique up to a Banach space
isomorphism. We shall later in this section prove a similar uniqueness result
for the operator space versions.

Our first result states:

\begin{theorem}
\label{thm1.1}

If $\sigma$ satisfies (\ref{eq1.1}) and (\ref{eq1.2}), then
$\tilde{X}_p(\sigma)^*$ is not Banach space isomorphic to a subspace of $S_{p'}$.
Consequently $\tilde{X}_p(\sigma)$ is not Banach space isomorphic to a complemented
subspace of $S_p$.
\end{theorem}

\bproof
Assume that $\tilde{X}_p(\sigma)^*$ is isomorphic to a subspace of $S_{p'}$ and
let $n\in \bN$ be given. By \cite[Corollary 8]{R} $\tilde{X}_p(\sigma)^*$
contains a basic sequence $(h_k)$ equivalent to the unit vector
basis of $\ell_2$ so that any $n$ elements of that sequence is
isometrically equivalent to the unit vector basis of $\ell_{p'}^n$.
From \cite[Proposition 4 and Lemma 1]{AL} it follows that $(h_k)$
has a subsequence which is 4-equivalent to the unit vector basis
of $\ell_2$. This is a contradiction for large $n \in \bN$.
\eproof

The next theorem is the operator space version of Rosenthal's lemma 7 in
\cite{R}.

\begin{proposition}
\label{thm1.1a}

Let $(g_n)$ be the natural basis of $X_{p*}(\sigma)$ and let $(E_j)$ be a
sequence of mutually disjoint finite subsets of $\bN$. For each $j \in
\bN$ we put
\begin{eqnarray}
f_j = \sum_{n\in E_j}\sigma_n ^{r/p}g_n \\
\beta_j = (\sum_{n \in E_j}\sigma_n ^r)^{\frac{1}{r}}\\
\tilde f_j = \beta_j ^{-r/p}f_j
\end{eqnarray}

$(\tilde f_j)$ is a cb-unconditional basic sequence, cb-isometrically
equivalent to
the natural basis of $X_{p*}(\beta)$ and there is a cb-contractive projection
of $X_{p*}$ onto $[f_j]$.
\end{proposition}

\bproof

We shall prove the theorem for $X_{p,c_p}(\sigma)$; the other cases can be proved
in a similar manner.

If $(x_j)_{j=1}^k \subseteq S_p$, then we get
\begin{equation*}
\|\sum_{j=1}^k x_j \otimes f_j \|_{S_p[X_{p,c_p}(\sigma)]} =
\|\sum_{j=1}^k \sum_{n \in E_j} \sigma_n ^{r/p}x_j \otimes [e_{nn}\oplus \sigma_n e_{n1}]\|_{S_p[X_{p,c_p}(\sigma)]}.
\end{equation*}

It easily follows that
\begin{equation*}
 \|\sum_{j=1}^k \sum_{n\in E_j} \sigma_n^{r/p}x_j \otimes e_{nn}\|_{S_p[S_p]}=\\ (\sum_{j=1}^k \|x_j\|^p
\sum_{n \in E_j}\sigma_n ^r)^{1/p} = \\
(\sum_{j=1}^k \|x_k \|^p \beta_j ^r)^{1/p}.
\end{equation*}

From (\ref{pc}) we get
\begin{equation*}
\|\sum_{j=1}^k \sum_{n \in E_j} \sigma_n^{r/p} x_j \otimes \sigma_n e_{n1} \|_{S_p[C_p]}=
\|(\sum_{j=1}^k \sum_{n \in E_j} \sigma_n ^{(\frac{2r}{p}+2)}x_j^*x_j)^{1/2}\|_{S_p} =
\|(\sum_{j=1}^k \beta_j^r x_j^*x_j)^{1/2}\|_{S_p}
\end{equation*}

and therefore
\begin{equation}
\label{eq1.3a}
\|\sum_{j=1}^k x_j \otimes \tilde f_j \|_{S_p[X_{p,c_p}(\sigma)]} =
\|\sum_{j=1}^k x_j \otimes [e_{jj} \oplus \beta_j e_{j1}]\|_{S_p[X_{p,c_p}(\beta)]}. 
\end{equation}

Together with Proposition \ref{prop01} this shows that $(\tilde f_j)$ is cb-isometrically equivalent to the natural basis $(g_j)$ of $X_{p,c_p}(\beta)$.

For all $x,y \in X_{p,c_p}(\sigma)$ we
put $<x,y>=\sum_{j=1}^{\infty} x(j)\overline{y(j)}\sigma_j^2$ (where $x(j)$, respectively y(j) denotes the j'th coordinate of $x$, respectively $y$ in the basis $(g_j)$)
and define
\begin{equation}
\label{eq1.3b}
Px = \sum_{j=1}^{\infty} <x,f_j>\beta^{-r} f_j \quad \quad \mbox{for
all $x \in X_{p,c_p}(\sigma)$}
\end{equation}

It follows immediately from Rosenthal's argumentation in \cite[Lemma 7]{R}
that in the Banach space sense $P$ is a contractive projection of
$X_{p,c_p}(\sigma)$ onto $[f_j]$. In addition we need to prove that $P$ is completely bounded with $\|P\|_{cb}=1$.

For every $n \in \bN$ we get
\begin{eqnarray}
\label{eq1.3c}
Pg_n = \sum_{j=1}^{\infty} <g_n,f_j>\beta_j^{-r}f_j & = &
\sigma_n^{r/p+2}\beta_{j_n}f_{j_n} = \\
\sigma_n^{r/p+2}\beta_{j_n}^{r/p-r}\tilde f_{j_n} & = &
\beta_{j_n}^{-r/p'}\sigma_n^{r/p'}\tilde f_{j_n}. \nonumber
\end{eqnarray}

where $j_n$ is chosen such that $n \in E_{j_n}$.

Let now $(x_n) \subseteq S_p$ be a finite sequence. From (\ref{eq1.3c})
and the first part of the proof we obtain
\begin{eqnarray}
& &  \|\sum_n x_n \otimes Pg_n \|_{S_p[X_{p,c_p}]} =
\|\sum_j \beta_j^{-r/p'}(\sum_{n \in E_j} \sigma_n^{r/p'}x_n)\otimes
\tilde f_j\|_{
S_p[X_{p,c_p}(\sigma)]}  \nonumber \\
& & = \| \sum_j \beta_j^{-r/p'}(\sum_{n \in
E_j}\sigma_n^{r/p'}x_n)\otimes [e_{jj}\oplus \beta_j
e_{j1}]\|_{S_p[X_{p,c_p}(\beta)]}. \label{eq1.3c2}
\end{eqnarray}

We estimate the two coordinates separately and start with:
\begin{eqnarray}
& & \|\sum_j \beta_j^{-r/p'}(\sum_{n \in E_j}
\sigma_n^{r/p'}x_n)\otimes e_{jj}\|_{S_p[S_p]} = (\sum_j
\beta_j^{-rp/p'} \|\sum_{n \in E_j}
\sigma_n^{r/p'}x_n \|_{S_p}^p)^{1/p} \label{eq1.3d} \\
& & \kl (\sum_j \beta_j^{-rp/p'}(\sum_{n \in E_j}\sigma_n^r)^{p/p'}
\sum_{n \in E_j} \|x_n \|_{S_p}^p)^{1/p} = (\sum_j \sum_{n \in
E_j}\|x_n\|_{S_p}^p)^{1/p} =  (\sum_n \|x_n\|_{S_p}^p)^{1/p}.
\nonumber
\end{eqnarray}

The estimate of the other coordinate is slightly more involved. For every
$\xi \in \ell_2$ and every $j$ we get
\begin{eqnarray*}
& & ((\sum_{n\in E_j}\sigma_n^{r/p'}x_n^*)(\sum_{n\in
E_j}\sigma_n^{r/p'}x_n)\xi,
\xi) =  \|\sum_{n\in E_j}\sigma_n^{r/p'}x_n \xi \|^2 \le (\sum_{n\in E_j}\sigma_n^{2r/p' -2})(\sum_{n\in E_j}\|\sigma_n x_n\xi\|^2)  \\
& & = \sum_{n\in E_j}\sigma_n^r \sum_{n\in E_j}(\sigma_n^2x_n^* x_n
\xi,\xi) =  \beta_j^r \sum_{n\in E_j}\sigma_n^2(x_n^*x_n \xi,\xi)
\end{eqnarray*}

which shows that in the sense of operators on $\ell_2$ we have:
\begin{equation*}
0 \le \sum_j \beta_j^{-r}(\sum_{n\in E_j}\sigma_n^{r/p'}x_n^*)
(\sum_{n\in E_j}\sigma_n^{r/p'}x_n) \le \sum_j \sigma_j^2x_j^* x_j.
\end{equation*}

Together with (\ref{pc}) and \cite[Theorem 2.3]{FH-GKP} this gives:
\begin{eqnarray}
& & \|\sum_j\beta_j^{-r/p'}\sum_{n\in E_j}\sigma_n^{r/p'}x_n \otimes
\beta_je_{j1}\|_{S_p[C_p]} =  \|(\sum_j\beta_j^{-r}(\sum_{n\in
E_j}\sigma_n^{r/p'}x_n^*) (\sum_{n\in
E_j}\sigma_n^{r/p'}x_n))^{1/2}\|_{S_p} \nonumber \\
& & = (tr([\sum_j
\beta_j^{-r}\sum_{n\in E_j}\sigma_n^{r/p'}x_n^*\sum_{n\in E_j}
\sigma_n^{r/p}x_n]^{p/2}))^{1/p}   \label{eq1.3e}  \\
& & \le  (tr([\sum_j \sigma_j^2x_j^*x_j]^{p/2}))^{1/p} =  \|\sum_j
x_j \otimes \sigma_j e_{j1}\|_{S_p[C_p]}. \nonumber
\end{eqnarray}

(\ref{eq1.3c2}), (\ref{eq1.3d}) and (\ref{eq1.3e}) show that $P$ is
completely bounded with $\|P\|_{cb} = 1$. \eproof

An application of Theorem \ref{thm1.1} shows like in the Banach space
case that if $\sigma$ in addition satisfies
(\ref{eq1.2}), then $X_{p*}(\sigma)$ is uniquely determined up to a
cb-isomorphism. This is the contents of the next theorem.

\begin{theorem}
\label{cor1.1a}
If $2< p < \infty$ and $\sigma$ and $\gamma$ are two sequences both satisfying
(\ref{eq1.1}) and (\ref{eq1.2}), then $X_{p*}(\sigma)$ is cb-isomorphic to
$X_{p*}(\gamma)$.
\end{theorem}

\bproof The proof follows the lines of the proofs of
\cite[Proposition 12 and Theorem 13]{R} and is based on
Pe{\l}czy\'{n}ski's decomposition method (see e.g. \cite[Theorem
2.a.3]{LT1}). We will therefore first prove that $X_{p*}(\gamma)$ is
cb-isomorphic to a cb-complemented subspace of $X_{p*}(\sigma)$ and
vice versa.

Since $\sigma$ satisfies (\ref{eq1.1}) and (\ref{eq1.2}), we can find a sequence
$(E_j)$ of mutually disjoint, finite subsets of $\bN$ so that
\begin{equation}
\label{eq1.21}
\gamma_j \le \beta_j = (\sum_{n\in E_j}\sigma_n^r)^{1/r} \le 2\gamma_j
\quad \quad \mbox{for all $j\in \bN$}
\end{equation}

From Proposition \ref{thm1.1a} it follows that $X_{p*}(\beta)$ is cb-isometric to
a subspace of $X_{p*}(\sigma)$ onto which there is a cb-contractive projection.
(\ref{eq1.21}) shows that $X_{p*}(\gamma)$ is 2-cb-isomorphic to
$X_{p*}(\beta)$.
By interchanging the roles of $\gamma$ and $\sigma$ we obtain that also
$X_{p*}(\sigma)$ is cb-isomorphic to a cb-complemented subspace of
$X_{p*}(\gamma)$.

The next step is to show that $X_{p*}(\sigma)$ is cb-isomorphic to
$X_{p*}(\sigma)\oplus X_{p*}(\sigma)$ but we shall only prove it for
$X_{p,c_p}(\sigma)$ since the other cases can be obtained in a similar manner.

(\ref{eq1.1}) and (\ref{eq1.2}) give that we can find a sequence
$\{E_{j,k} \mid j\in \bN, \quad  k\in \bN \}$ of mutually disjoint finite subsets of $\bN$ so that
\begin{equation}
\label{eq1.22}
\sigma_j \le \beta_{j,k} = (\sum_{n\in E_{j,k}} \sigma_n^r)^{1/r} \le
2\sigma_j \quad \quad \mbox{for all $j,k \in \bN$}
\end{equation}

Put $\beta_k =(\beta_{j,k})_{j=1}^{\infty}$, let
$\tilde f_{j,k}= \beta_{j,k}^{-r/p}\sum_{n\in E_{j,k}}\sigma_n^{r/p}e_{nn}\oplus
\sigma_n e_{n1}$ and define \mbox{$Z= [\tilde f_{j,k} \mid j,k \in \bN]$},
\mbox{$Z_1 = [\tilde f_{j,k} \mid j\in \bN \quad k \ge 2]$}. It follows
from Proposition \ref{thm1.1a} that $Z$ is cb-contractively complemented in
$X_{p,c_p}(\sigma)$ and that for all $k \in \bN$ $[\tilde f_{j,k}]$ is cb-contractively complemented and cb-isometric to $X_{p,c_p}(({\beta}_k))$ which in turn
is 2-cb-isomorphic to $X_{p,c_p}(\sigma)$. Hence $Z$ can be viewed as an
infinite direct sum of copies of $X_{p,c_p}(\sigma)$.
Let $T: span \{\tilde f_{j,k} \mid j,k \in \bN \} \rightarrow Z_1$ be defined by $T \tilde f_{j,k} = \tilde f_{j,k+1}$. We shall show that $T$ extends to
a cb-isomorphism of $Z$ onto $Z_1$. If $(x_{j,k}) \subseteq S_p$ is a finite
sequence, then we get from (\ref{eq1.22}) and \cite[Theorem 2.3]{FH-GKP} that
\begin{eqnarray}
& & \|(\sum_k \sum_j \beta_{j,k+1}^2 x_{j,k}^* x_{j,k})^{1/2}
\|_{S_p} \kl
2\|(\sum_k \sum_j \sigma_j^2 x_{j,k}^* x_{j,k})^{1/2}\|_{S_p}  \label{eq1.18} \\
& & \le 2\|(\sum_k \sum_j \beta_{j,k}^2 x_{j,k}^*
x_{j,k})^{1/2}\|_{S_p}.  \nonumber
\end{eqnarray}

In the same manner we get
\begin{equation}
\|(\sum_k \sum_j \beta_{j,k}^2 x_{j,k}^* x_{j,k})^{1/2}\|_{S_p} \le
2\|(\sum_k \sum_j \beta_{j,k+1} x_{j,k}^* x_{j,k})^{1/2}\|.
\end{equation}

Similar estimates can easily be obtained for the corresponding p-norms which implies that we have
\begin{equation*}
\frac{1}{2}\|\sum_k \sum_j x_{j,k} \otimes \tilde f_{j,k}\|_{S_p[X_{p,c_p}(\sigma)]}
\le \|\sum_k \sum_j x_{j,k} \otimes \tilde f_{j,k+1}\|_{S_p[X_{p,c_p}(\sigma)]} \le
2\|\sum_k \sum_j x_{j,k} \otimes \tilde f_{j,k}\|_{S_p[X_{p,c_p}(\sigma)]},
\end{equation*}

which shows that $T$ can be extended to a cb-isomorphism of $Z$ onto $Z_1$.

Letting $\sim_{cb}$ denote ``cb-isomorphic to'', we obtain from the above that
$Z \sim_{cb} X_{p,c_p}(\sigma) \oplus Z$. Since $Z$ is cb-complemented in
$X_{p,c_p}(\sigma)$, we can find a closed subspace $U \subseteq X_{p,c_p}(
\sigma)$ so that $X_{p,c_p}(\sigma)= Z \oplus U \sim_{cb} X_{p,c_p}(\sigma)
\oplus Z \oplus U \sim_{cb} X_{p,c_p}(\sigma) \oplus X_{p,c_p}(\sigma)$.

We are now ready to show that $X_{p,c_p}(\gamma)$ is cb-isomorphic to
$X_{p,c_p}(\sigma)$. Indeed, since by the above $X_{p,c_p}(\gamma)$ is
cb-isomorphic to a cb-complemented subspace of $X_{p,c_p}(\sigma)$, we can
find a closed subspace $G \subseteq X_{p,c_p}(\sigma)$ so that
\begin{equation*}
X_{p,c_p}(\sigma) \sim_{cb} X_{p,c_p}(\gamma)\oplus G \sim_{cb}
X_{p,c_p}(\gamma) \oplus X_{p,c_p}(\gamma)\oplus G \sim_{cb}
X_{p,c_p}(\gamma)\oplus X_{p,c_p}(\sigma) \sim_{cb} X_{p,c_p}(\gamma)
\end{equation*}

where the last $\sim_{cb}$ follows by interchanging the roles of $\sigma$ and
$\gamma$.

\eproof

Exploiting the decomposition method a bit more we can actually obtain that
also the space $Z$ in the above proof is cb-isomorphic to
$X_{p,c_p}(\sigma)$.

We are now going to define some operator spaces which we shall call
matricial Rosenthal spaces.

The space $\tilde{Y}_p(\sigma)$ is defined to be the subspace of ${\K}_p
\oplus_p (\sum_{n=1}^{\infty} S_2 ^n)_2$ consisting of all elements of the
form $((x_n,\sigma_n x_n))$ where $x_n \in S_{p}^n$ for all $n\in \bN$,
i.e we require:
\begin{equation}
\label{eq1.5}
\sum_{n=1}^{\infty} \|x_n\|_{S_p ^n}^p < \infty
\quad \mbox{and} \quad \sum_{n=1}^{\infty} \sigma_n ^2 \|x_n\|_{S_2
^n}^2 < \infty.
\end{equation}

We can view $(\sum_{n=1}^{\infty} S_2 ^n)_2$ isometrically as a subspace
of $C_p[C_p]$ in the following way: Choose a sequence $(m_n)$ of
integers so that $m_1=0$ and $m_{n+1}-m_n =n$ for all $n \in \bN$. If 
$x=(x_n) \in (\sum_{n=1}^{\infty} S_2 ^n)_2$ with $x_n =
(t_{ij}^n)_{i,j=1}^n$, we can identify $x$  with 
$\sum_{n=1}^{\infty}\sum_{i,j=m_n+1}^{m_{n+1}} t_{ij}^ne_{ij} \in
C_p[C_p]$. Similarly we can consider $(\sum_{n=1}^{\infty} S_2 ^n)_2$ as
a subspace of $R_p[R_p]$, respectively of $C_p[C_p] \oplus_p R_p[R_p]$. 

Hence there is a canonical Banach space isometry $w_{\sigma}$ of 
$\tilde{Y}_p(\sigma)$ into the operator space 
${\K}_p \oplus_p C_p[C_p]$ and we put
$Y_{p,c_p} = w_{\sigma}(\tilde{Y}_p(\sigma))$. Similarly we define the
spaces $Y_{p,r_p}(\sigma)$ and $Y_{p,c_p \cap r_p}(\sigma)$. In the rest
of this paper we shall put $Y_p(\sigma)=Y_{p,c_p \cap r_p}(\sigma)$. 

In the sequel we often have to consider cb-maps to or from these spaces
and it is therefore worthwhile to mention how the norm in
$S_p[Y_{p,c_p}(\sigma)]$ is computed (the other cases follow similarly). Let us
just compute the ``column part'' of $S_p[Y_{p,c_p}(\sigma)]$. To this
end let $X_n \in S_p \otimes S_p^n$ for all $n \in \bN$. We can then
find $(x_{jk}^n)\in S_p^n$ so that 
\begin{equation*}
X_n = \sum_{j,k= m_n+1}^{m_{n+1}} x_{jk}^n \otimes e_{jk}
\end{equation*}
for every $n \in \bN$. Note that
\begin{equation}
\label{egyp0}
X_n^*X_n = \sum_{k,l=m_n+1}^{m_{n+1}}(\sum_{j=m_n+1}^{m_{n+1}} x_{jk}^{n*}
x_{jl}^n)e_{kl}.
\end{equation}
Using Proposition \ref{prop02} we get that:
\begin{eqnarray}
\label{eqyp}
\|\sum_n\sigma_n X_n\|_{S_p[C_p[C_p]]} & = & \|\sum_n
\sum_{j,k=m_n+1}^{m_{n+1}} x_{jk}^n\otimes e_{jk}\|_{S_p[C_p[C_p]]}\\
= \|(\sum_n\sigma_n^2\sum_{j,k=m_n+1}^{m_{n+1}}x_{jk}^{n*}x_{jk}^n)^{\frac{1}{2}}\|_{S_p} & = & 
\|(\sum_n \sigma_n^2 (id\otimes tr)(X_n^*X_n)^{\frac{1}{2}}\|_{S_p},
\nonumber
\end{eqnarray} 

where we have used (\ref{egyp0}) to get the last equality.
Comparing this with the similar calculations for the other cases it is
readily verified that $Y_p(\sigma)$, $Y_{p,c_p}(\sigma)$, and $Y_{p,r_p}(\sigma)$
are mutually non-cb-isomorphic.

Since ${\K}_p \oplus_p C_p[C_p]$ is cb-isomorphic
to a subspace of $S_p$ the same holds for $Y_{p,c_p}(\sigma)$ as
well. In a similar manner we get that $Y_{p,r_p}(\sigma)$ and
$Y_p(\sigma)$ are cb-isomorphic to subspaces of $S_p$.
We have the following result on these spaces.

\begin{theorem}
\label{thm1.2}
Both ${\K}_p$ and $X_{p,c_p}(\sigma)$ (respectively $X_{p,r_p}(\sigma)$) are
cb-isomorphic to complemented subspaces of $Y_{p,c_p}(\sigma)$
(respectively $Y_{p,r_p}(\sigma)$).  Consequently $\tilde{Y}_p(\sigma)$
is not Banach space isomorphic to a complemented subspace of $S_p$ if $\sigma$
satisfies (\ref{eq1.1}) and (\ref{eq1.2}).
\end{theorem}

\bproof Let $U=X_{p,c_p}(\sigma)$ (respectively
$U=X_{p,r_p}(\sigma)$) and $W=Y_{p,c_p}(\sigma)$ (respectively
$W=Y_{p,r_p}(\sigma)$). If $(n_k)\subseteq \bN$ is a sequence so
that $\sum_{k=1}^{\infty} \sigma_{n_k}^{\frac{2p}{p-2}} < \infty$,
then the subspace $V$ consisting of those $(x_n,\sigma_n x_n)\in W$
for which $x_n = 0$ for all $n\neq n_k$ is readily seen to be
completely complemented by a projection of cb-norm one and
completely isomorphic to ${\K}_p$.

It is obvious that $U$ can be identified cb-isometrically
with the subspace of $W$ consisting of those $(x_n,\sigma_n x_n) \in \tilde{Y}_p(\sigma)$ for
which $x_n$ is a one-dimensional operator on $\ell_2$ for all $n\in \bN$. This space is
clearly the range of a cb-contractive projection.

It now follows directly from Theorem \ref{thm1.1} that $\tilde{Y}_p(\sigma)$
cannot be Banach space isomorphic to a complemented subspace of $S_p$ if
$\sigma$ satisfies (\ref{eq1.1}) and (\ref{eq1.2}).
\eproof

The last spaces we are going to investigate are defined as follows:

\begin{equation}
\label{eq1.6}
Z_{p,c_p}(\sigma) = \{(x,xD_{\sigma}) \mid x\in A_{\sigma} \} \subseteq
S_p \oplus_p C_p[C_p].
\end{equation}

\begin{equation}
\label{eq1.6a}
Z_{p,r_p}(\sigma) = \{(x,D_{\sigma}x) \mid x\in A_{\sigma} \} \subseteq
S_p \oplus_p R_p[R_p].
\end{equation}

\begin{equation}
\label{eq1.6b}
Z_p(\sigma) = \{(x,xD_{\sigma},D_{\sigma}x) \mid x\in A_{\sigma} \} \subseteq
S_p \oplus_p C_p[C_p ]\oplus_p R_p[R_p].
\end{equation}
where
\[
A_{\sigma} = \{x\in S_p \mid xD_{\sigma} \in S_2 \}.
\]

In (\ref{eq1.6}) we consider $xD_{\sigma}$ as an element of $C_p[C_p]$
and similarly in (\ref{eq1.6a}) and (\ref{eq1.6b}).

In the sequel we let $Z_{p,*}(\sigma)$ denote any of these
spaces. Clearly they are isomorphic as Banach spaces, are mutually
non-cb-isomorphic and cb-embed into $S_p$.

The next theorem gives the basic properties of the spaces $Z_{p*}(\sigma)$.

\begin{theorem}
\label{thm1.3}
The space $Z_{p,*}(\sigma)$ has the following properties:
\begin{itemize}

\item[(i)] If $\sigma$ satisfies (\ref{eq1.1}), then $S_p$ is cb-isomorphic
to a cb-complemented subspace of
$Z_{p,*}(\sigma)$.

\item[(ii)] If $\sigma$ satisfies both (\ref{eq1.1}) and (\ref{eq1.2}), then
$Z_{p,*}(\sigma)$ is not isomorphic to a complemented
subspace of $S_p$.
\end{itemize}
\end{theorem}

\bproof
(i): We shall only give the argument for $Z_{p,c_p}(\sigma)$. The proof
for other spaces can be made in a similar manner. Let $(n_k) \subseteq
\bN$ be a sequence so that $\sum_{k=1}^{\infty} \sigma_{n_k}^{\frac{2p}{p-2}} < \infty$ and
let $V$ consist of those $(x,xD_{\sigma})\in Z_p(\sigma)$ for
which $x_{ij} =0$ unless $j=n_k$ for some $k\in \bN$. It is readily
verified that $V$ is cb-isomorphic to $S_p$. From Arazy
\cite[Theorem 1.1]{A} it follows that $V$ contains another subspace
$U$ cb-isomorphic to $S_p$ and which is complemented in $Z_p(\sigma)$.
This shows (i).

(ii): $X_p(\sigma)$ can easily be identified with those
$(x,xD_{\sigma})\in Z_p(\sigma)$ for which $x$ is a diagonal matrix.
This subspace is clearly contractively complemented in
$Z_p(\sigma)$. It now follows from Theorem \ref{thm1.1} that
$Z_p(\sigma)$ is not isomorphic to complemented subspace
of $S_p$.
\eproof

Before we go on we need the following lemma on non-commutative
$L_p$-spaces.

\begin{lemma}
\label{lemma1.4}
Let $1<p< \infty$ and let $\cN$ be a von Neumann algebra so that
$L_p(\cN)$ is separable and $L_p(0,1)$ does not embed
isomorphically into $L_p(\cN)$. Then there exist sequences $(I_k)$
of countable sets and $(n_k)\subseteq \bN$ so that
\begin{equation}
\label{eq1.7a}
L_p(\cN) = (\sum_{k=1}^{\infty} \ell_p(I_k,S_p^{n_k}))_p.
\end{equation}
\end{lemma}

\bproof
Since $L_p(0,1)$ does not embed into $L_p(\cN)$, it follows from a
a result of Marcolino \cite{Ma} that $\cN$ is a type $I$ factor and
therefore the separability of $L_p(\cN)$ and \cite{Ta} give that there exist
measure spaces $(\Omega_k,\Sigma_k,\mu_k)$
for all $k \in \bN$ and $(n_k)\subseteq \bN$ so that
\begin{equation}
\label{eq1.7}
L_p(\cN) = (\sum_{n=1}^{\infty} L_p(\Omega_k,\Sigma_k,\mu_k,S_p^{n_k})_p.
\end{equation}
Again, since $L_p(0,1)$ does not embed into $L_p(\cN)$, it follows that
all the measure spaces on the right side of \eqref{eq1.7} are purely
atomic.
\eproof

We are now able to prove:
\begin{theorem}
\label{thm1.5}
If $\sigma$ satisfies (\ref{eq1.1}) and (\ref{eq1.2}), then none of the
spaces $X_p(\sigma)$, $Y_p(\sigma)$ or $Z_p(\sigma)$
are isomorphic to an $L_p(\cN)$-space where $\cN$ is a von Neumann
algebra.
\end{theorem}

\bproof
Let $V$ be one of the spaces above and assume that there exists
von Neumann algebra $\cN$ so that $V$ is isomorphic to $L_p(\cN)$. Since
it follows from \cite[Theorem 6]{AL} that $L_p(0,1)$ does not
embed into $S_p$, $L_p(\N)$ has the form of \eqref{eq1.7a} by
Lemma \ref{lemma1.4} and therefore it is isomorphic to a
complemented subspace of $S_p$. This contradicts Theorems
\ref{thm1.1}, \ref{thm1.2} and \ref{thm1.3} above.
\eproof

%% file: njnpart2new11.tex
\section{The operator space structure of the  classical Rosenthal sequence spaces}
\label{sec2}
\setcounter{equation}{0}

In this section we wish to discuss the operator space
structure of the Rosenthal sequence spaces defined in Section 1 and it
turns out that the local structure of these spaces behaves quite differently.
However, due to the non-commutative Burkholder-Rosenthal inequalities
 \cite{JX}, \cite{JXII} the probabilistic viewpoint from the commutative case is still adequate to determine this
structure.

Let $(\si_{i})$ be a sequence such that $0\le \si_i \le
1$ and let $A_i\subset [0,1]$, $i\in \bN$ be intervals of
measure $\mu(A_i)=\si_i^r$, where
$\frac12=\frac1p+\frac1r$. We define
$f_i((t_j))= \mu(A_i)^{-\frac1p}
1_{A_i}(t_i)$ for all sequences $(t_j) \subseteq [0,1]$. The sequence $(f_i)_{i\in \nz}$
is a sequence of independent random variables on $[0,1]^{\bN}$. For
sequences $(s_i)$ with  finite support we define
 \for
 u((s_i)) &=& \summ_{i=1}^\8 s_i \eps_i f_i \pl ,\\
 u_c((s_i)) &=& \summ_{i=1}^\8 s_i e_{i,1} \pl \eps_i f_i
 \pl ,\\
 u_r((s_i)) &=& \summ_{i=1}^\8 s_i e_{1,i}\pl \eps_i f_i  \pl ,
 \mel
where $(\varepsilon_i)$ denotes the sequence of Rademacher functions on
$[0,1]$.

Following Rosenthal's argument from \cite{R} using \cite{JXII} we can now
obtain

\begin{proposition} \label{brr} Let $2\le p<\8$, then $u$, $u_c$, $u_r$
is a cb-isomorphism between $X_p(\si)$,
$X_{p,c_p}(\si)$ and $X_{p,r_p}(\si)$ and the image
of $u$ in $L_p([0,1]^{\bN}$, $u_c$ in $L_p([0,1]^{\bN};C_p)$,
$u_r$ in $L_p([0,1]^{\bN};R_p)$, respectively. The images are cb-complemented in the
respective spaces.
\end{proposition}

\bproof
We shall only prove the proposition for $u_c$ since the other cases go
similarly. Let $(x_i)_{i=1}^n \subseteq S_p$ be arbitrary. From
\cite[Corollary 1.5]{JXII} and Proposition \ref{prop02} we get letting $\sim$ denote
two-sided inequalities with constants only depending on $p$:

\begin{eqnarray}
& &\|\sum_{i=1}^n x_i \otimes \varepsilon_i f_i e_{i1}\|_{S_p[L_p((0,1);C_p)]} \nonumber \\
& & \sim \max\{(\sum_{i=1}^n \|x_i\|_{S_p}^p \|f_i\|_p
^p)^{\frac{1}{p}}, \|(\sum_{i=1}^n x_i^*x_i {\mathbb
E}(f_i^2))^{\frac12}\|_{S_p}, (\sum_{i=1}^n \|x_i\|_{S_p}^p{\mathbb
E}(f_i^2)^{\frac p2})^{\frac{1}{p}}\}
  \label{njneq2.1} \\
& & \sim  \|\sum_{i=1}^n x_i \otimes (e_{ii} \oplus \sigma_i
e_{i1})\|_{S_p[X_{p,c_p}(\sigma)]}
 \nonumber
\end{eqnarray}

where we in the last equivalence have used that for all $1\le i \le n$ we
have $\|f_i\|_p =1$, ${\mathbb E}(f_i^2)= \sigma_i^2$ and
${\mathbb E}(f_i^2)^{\frac{p}{2}}= \mu(A_i)^{\frac{p}{2}-1}\le 1$. By
Lemma \ref{prop01} $u_c$ is a cb-isomorphism.

For every $1\le i \le n$ we put $f'_i = \mu(A_i)^{\frac{1}{p'}}1_{A_i}$ and
$u_{p'}((s_i))= \sum s_i \varepsilon_i f'_i$. Using the second part of
\cite[Theorem 0.1]{JXII} in a similar manner as above we achieve that
$u_{p'}$ acts as a cb-bounded operator from $X_{p,c_p}^*$ to $L_{p'}(0,1)$.
It is readily verified that $u_c u_{p'}^*$ is a cb projection of $L_p(0,1)$
onto the range of $u_c$.
\eproof

\begin{cor} 
\label{njn2.2}
The space $X_p(\si)$, $X_{p,c_p}(\si)$ and
$X_{p,r_p}$ have the $\gamma_p$-AP. More precisely, $X_p(\si)$
admits an approximate diagram
    \[ \begin{array}{ccccc}  X_p& &\stackrel{id}{\longrightarrow}& &X_p  \\
                                & {\scriptstyle v_n}\searrow & & \nearrow {\scriptstyle w_n}& \\
                               & & \ell_p^{n_k}& &
        \end{array} \]
For $X_{p,c_p}(\si)$ and $X_{p,r_p}(\si)$ we have to replace $\ell_p^{n_k}$ by
$\ell_p^{n_k}(C_p^{n_k})$ and $\ell_p^{n_k}(R_p^{n_k})$, respectively.
\end{cor}

\begin{cor} If $\sigma$ satisfies (\ref{eq1.1}), then the Rosenthal
spaces $X_p(\si)$ are ${\cC\cO\cL}_p$-spaces.
\end{cor}

\bproof
Follow the proof of \cite[Proposition 2.4]{JNRX}, using Corollary \ref{njn2.2}
and the fact that $X_p(\sigma)$ contains completely complemented copies of 
$\ell_p^n$'s far out.
\eproof

In the following we want to show that the
Rosenthal spaces $X_{p,c_p}(\si)$ and
$X_{p,r_p}(\si)$ are no longer ${\cO\cL}_p$.
Indeed, the mixture between the Hilbert space
structure and the $\ell_p$ structure forms the
crucial obstacle.

\begin{lemma}\label{cpp}
If $1\le p< \8$ and $\N$ is a finite von Neumann
algebra, then $C_p$ is not cb-isomorphic to a
subspace of $R_p(L_p(\N))$. Similarly, $R_p$ is not
cb-isomorphic to a subspace of $C_p(L_p(\N))$.
\end{lemma}

{\bf Proof:} Assume to the contrary that $C_p$ is isomorphic
to a subspace of $R_p(L_p(\N))$. Using the natural
isomorphism between  $R_p(R_p)$ and $R_p$, we
deduce that $S_p=R_p(C_p)$ is a Banach space
isomorphic to a subspace of $R_p(L_p(\N)) \subset
L_p(B(\ell_2)\ten \N)$. However, for $x\in
R_p(L_p(\N))$ and $p\ge 2$, we have
 \[ \noo x\rrm_2 \lel \noo xx^*\rrm_{L_1(\N)}^{\frac12} \kl
 \noo xx^*\rrm_{\frac{p}{2}}^{\frac12} \kl \noo
 x\rrm_p \pl .\]
Thus $R_p(L_p(\N))$ is isomorphic to a subspace of
$L_p(B(\ell_2)\ten \N)\cap L_2(B(\ell_2) \ten \N\
)$ for $2\le p<\8$. For $1\le p\le 2$ a similar
argument shows  that $R_p(L_p(\N))$ is isomorphic
to a subspace of $L_p(\N\ten B(\ell_2)) + L_2(\N\ten
B(\ell_2))$. According to \cite{J2} these spaces
are isomorphic to complemented subspaces of
$L_p(\M)$ for some finite von Neumann algebra $\M$.
Hence, $S_p$ is isomorphic to a subspace of
$L_p(\M)$. This contradicts Suckochev's result for
$p\ge 2$, \cite{Su}, or \cite{HRS} for $1\le p\le
2$. By symmetry the same holds for $R_p$ and
$C_p$ interchanged. \qed

\begin{cor} \label{finit}
Let $2<p,r<\8$ and $\frac12=\frac1p+\frac 1r$. If
$\si \notin \ell_r$, then the spaces
$X_{p,c_p}(\si)$ and $X_{p,r_p}(\si)$ are not
cb-isomorphic to subspaces of $L_p(\N)$ with $\N$
finite.
\end{cor}

{\bf Proof:} Assume first that there is an infinite set  $A\subset \nz$
so that $\inf_{k\in A} \si_k>0$. By interpolation we
deduce  that for the bounded sequence
$(\si_k^{-1})_{k\in A}$ the diagonal map
$D_{\si^{-1}}:C_p\to \ell_p$  is completely bounded.
Hence, the subspace of $X_{p,c_p}(\si)$ consisting of the
sequences having their support in $A$ is cb-isomorphic
to $C_p$. In particular it cannot embed into
$L_p(\N)$ cb-isomorphically. Thus $X_{p,c_p}(\si)$
can not embed either in this case. Since $\sum_j
\si_j^r=\8$, we can in the general case find disjoint finite subsets
$A_j$ such that if
 \[ \beta_j  \lel \kla \sum_{i\in A_j} \si_i^r \mer^{\frac1r} \pl, \]
then $\inf\beta_j > 0$. Proposition \ref{thm1.1a} gives that
$X_{p,c_p}(\beta)$ is cb-isomorphic to a subspace of
$X_{p,c_p}(\si)$ and by the above cb-isomorphic to $C_p$ and hence
the assertion follows. A similar argument applies for
the  row spaces.\qed

\begin{lemma}\label{cup}
If $1\le p\le \8$, then $\prod_\U \ell_p$ is
completely isometrically isomorphic to $L_p(\N)$
for a commutative von Neumann algebra $\N$.
\end{lemma}

{\bf Proof:} Let $\N=(\prod_\U \ell_1)^*$. According to
Raynaud's Theorem \cite{Re} we deduce that for all $n\in \nz$
$(S_1^n(\prod_U \ell_1))^*\lel M_n(\N)$ where $\N$ is a commutative
von Neumann algebra obtained as the weak closure of $\prod \ell_\8$.
Together with \cite[Lemma 5.4]{P1} this
implies that

  \[ L_p(M_n \ten \N ) \lel \prod S_p^n(\ell_p) \lel
  S_p^n(\prodd_\U \ell_p) \lel S_p^n(L_p(\N))  \pl .\]
Thus $\prod_U L_p$ is completely isometrically
isomorphic to $\ell_p(\N)$.\qed

Our aim is now to show that $X_{p,c_p}(\sigma)$ is not a
rectangular ${\cO\cL}_p$-space.

\begin{lemma} \label{dist}
If $2\le p\le \8$, then for all $\nen$
 \[ n^{\frac12-\frac1p} \kl   \inf_{E\subset C_p(L_p(0,1))} d_{cb}(R_p^n,E) \kl c_p \pl  n^{\frac12-\frac1p} \pl .\]
The same estimates hold if $R_p$ and  $C_p$
are interchanged.
\end{lemma}

{\bf Proof:} By interpolation
 \[ d_{cb}(R_p^n,R_p^n\cap C_p^n) \kl \noo id:R_p^n\to C_p^n\rrm_{cb}
 \noo id:R_p^n\cap C_p^n\to R_p^n\rrm_{cb} \kl n^{\frac{1}{2}-\frac{1}{p}} \pl .\]
By the non commutative Khintchine inequality \cite{Lu}
 \[  d_{cb}(R_p^n\cap C_p^n,{\rm span}\{g_j| j=1,..,n\}) \kl c_p \pl ,\]
where the $g_j$'s are independent Gaussian variables. To prove the lower
estimate, we consider  $E\subset
L_p(C_p)$ and a complete contraction $\phi:R_p^n\to
E$ and an isomorphism.  Let $x_i=\phi(e_{1i})$, then
 \for
 & & \kla   \intt \kla \summ_{i=1}^n \noo x_i(s)\rrm_2^2 \mer^{\frac{p}{2}} d\mu(s) \mer^{\frac1p}
 \lel   \noo \summ_{i=1}^n e_{i,1}\ten
 x_i\rrm_{L_p(C_p^n(C_p))} \\
 & & \kl     \noo \phi\rrm_{cb}\noo \summ_{i=1}^n e_{i,1}\ten
  e_{1,i}\rrm_{C_p^n[R_p^n]} \lel
  \noo id\rrm_{S_p^n} \lel n^{\frac1p} \pl .
 \mel
However, this implies
 \for
 \sqrt{n} &=& \kla \ez  \noo \summ_{i=1}^n \eps_i e_{1,i}\rrm_2^2 \mer^{\frac12}
  \lel  \kla \ez  \noo \summ_{i=1}^n \eps_i \phi^{-1}(x_i)\rrm_2^2 \mer^{\frac12} \\
 &\le& \noo \phi^{-1} \rrm \pl
 \kla \ez  \noo \summ_{i=1}^n \eps_i (x_i)\rrm_{L_p(\ell_2)}^2 \mer^{\frac12} \\
 &\le& \noo \phi^{-1} \rrm \pl  \kla \intt (\ez \noo \summ_{i=1}^n \eps_i x_i(s)\rrm_2^2)^{\frac p2} \mu(s) \mer^{\frac1p} \\
 &=& \noo \phi^{-1} \rrm \pl   \kla \intt \kla \summ_{i=1}^n \noo x_i(s)\rrm_2^2\mer^{\frac p2} \mu(s) \mer^{\frac1p}
 \kl
  \noo \phi^{-1} \rrm \pl n^{\frac1p} \pl .
 \mel
The assertion is proved.\qed

Using a similar idea we can even prove a slighly
stronger statement

\begin{lemma} \label{dist2}
If $2\le p\le \8$, then for all $\nen$
 \[  \frac{1}{c_p} \pl n^{\frac12-\frac1p} \kl   \inf_{E\in QS(\prod_\U L_p(C_p))} d_{cb}(R_p^n,E) \kl c_p \pl  n^{\frac12-\frac1p} \pl .\]
Here $c_p$ is an absolute constant and  $QS(\prod_\U
L_p(C_p))$ stands for the class of quotients of
subspace of ultraproducts of $C_p(L_p(0,1))$. The
same estimates holds exchanging  $R_p$ with $C_p$.
\end{lemma}

{\bf Proof:}  Let $T:C_p^n\to L_p(0,1)$ be defined
by $T(e_{i1})=\eps_i$, where $(\eps_i)_{i=1}^n$ are
Bernoulli random variables.  We will use a  a
sequence of independent normalized complex gaussian
random variables $(g_j)$ on $(\Om',\mu')$.  Let
$h_1,..,h_n\in L_p(\Om,\mu;\ell_2)$. Then, we
deduce from the Khinchine/Kahane's inequality
\cite{Kah}
 \for
 \noo \summ_{i=1}^n \eps_i h_i\rrm_{L_p(\ell_2)}
 &=& \noo g_1\rrm_p^{-1} \pl \kla \intt_{\Om\times \Om'} \intt_0^1 \bet
  \summ_{i=1}^n \summ_{=1}^\8   \eps_i(s) g_j(\om') h_i(j,\om))\rag^p \pl ds d\mu'(\om') d\mu(\om) \mer^{\frac1p} \\
 &\le&  \noo g_1\rrm_p^{-1} \pl c_0 \sqrt{p} \pl
  \kla \intt_{\Om\times \Om'} \kla \summ_{i=1}^n
 \bet \summ_{j=1}^\8  g_j(\om') h_i(j,\om))\rag^2\mer^{\frac{p}{2}}
    d\mu'(\om') d\mu(\om) \mer^{\frac1p} \\
 &\le&  \noo g_1\rrm_p^{-1} \pl c_0^2 p \pl
  \kla \intt_{\Om} \kla \summ_{i=1}^n \summ_{j=1}^\8
  \bet  h_i(j,\om))\rag^2\mer^{\frac{p}{2}}
 d\mu(\om) \mer^{\frac1p}  \pl .
 \mel
Since for $p\ge 2$, we have $\noo g_1\rrm_p \sim \sqrt{p}$ we
deduce
 \[  \noo  T\ten  id_{C_p(L_p(\Om))}: C_p^n(C_p(L_p(\Om))) \to C_p(L_p([0,1]\times \Om)) \rrm \kl c_0^3 \sqrt{p} \pl .\]
This remains true if we pass to an ultraproduct and
then to  a quotient of a subspace.   On the other
hand, we have seen in Lemma \ref{dist} that
 \[ \noo T\ten id_{R_p^n} \rrm \gl n^{\frac12-\frac1p} \pl .\]
Therefore the distance is bigger that
$\frac{n^{\frac12-\frac1p}}{c_0^3\sqrt{p}}$.\qed

The next lemma is a kind of ``folklore'' but for the convenience of the
reader we give a proof.

\begin{lemma} \label{leftm}
 Let $\M$ be a von Neumann algebra
and $2< p\le \8$, $2\le r<\8$ such that
$\frac12=\frac1p+\frac1r$. Let $F\subset L_p(\M)$
be a subspace and  $T:F\to R_p$ be a linear map.
$T$ is a complete contraction if and only if there
exists a norm one elements $a\in L_r(\M)$ and a
contraction $W:L_2(\M)\to \ell_2$ such that
 \[ T(x)\lel W(ax) \]
for all $x\in L_p(\N)$. In particular, $T$ admits a
completely contractive extension
$\hat{T}:L_p(\M)\to R_p$. Similarly, every complete
contraction $T:F\to C_p$ has a completely
contractive extension of the form $T(x)=W(xa)$.
\end{lemma}

{\bf Proof:} Let $(x_j)$ be a finite sequence in
$F$, then
 \for
 \kla \summ_j \noo T(x_j)\rrm_2^2 \mer^{\frac12}
& = & \noo \summ_j e_{j,1}\ten T(x_j)\rrm_{R_p(R_p)}
 \le  \noo \summ_j e_{j,1} \ten x_j\rrm_{R_p(L_p(\M))} \\
 &=&  \noo \summ_j x_jx_j^*\rrm_{\frac{p}{2}}
 \lel  \sup_{a\ge 0, \noo a\rrm_{\frac{r}{2}} \le 1} \kla \summ_j tr(ax_jx_j^*)\mer^{\frac12} \pl .
 \mel
Let $B$ be positive part of the unit ball of $L_{\frac r2}(\M)$.
The function $f_x(a)\mapsto tr(ax^*x)$ is continuous with
respect to the weak$^*$ topology. Hence, the standard separation
yields a probability measure $\mu$ on $B$ such that
 \[ \noo T(x)\rrm_2^2 \kl \intt_B tr(ax^*x) d\mu(a) \lel tr((\intt_B a d\mu(a))  x^*x) \pl .\]
By convexity, $b=(\intt_B a d\mu(a))\in B$ and
therefore
 \[ \noo T(x)\rrm_2 \kl \noo b^{\frac12}x\rrm_2 \pl .\]
Let $H=\{b^{\frac12}x\p|\p x\in F\}\subset
L_2(\M)$. Thus there is a linear contraction $W_1:H
\to \ell_2$ such that $W_1(b^{\frac12}x)=T(x)$. If
$P$ denotes the orthogonal projection onto $H$, then
$W=W_1P$ satisfies the assertion. To prove the
converse, we assume $T(x)=W(ax)$ for some $a\in
L_r(\M)$ of norm less than  one. Let
$L_a:L_p(\M)\to L_2(\M)^{r_p}$ be the left
multiplication $L_a(x)=ax$. Let $\phi: L_{\frac
p2}\to \cz$ be the induced linear functional
$\phi(y)=tr(ya^*a)$ of norm less than one. If
$x\in L_p(B(\ell_2)\ten \M)$, we deduce that
for every functional the $cb$-norm coincides with
the norm
 \for
  & & \noo (id \ten L_a)(x)\rrm_{S_p(L_2(\M)^{r_p})}
  \lel  \noo (id \ten tr)((a\ten id)x x^*(a^*\ten id)) \rrm_{S_{\frac{p}{2}}}^{\frac12} \\
  & & \lel  \noo (id \ten tr)(x x^*(a^*a\ten id)) \rrm_{S_{\frac{p}{2}}}^{\frac12}
    \lel   \noo (id \ten \phi)(xx^*)\rrm_{S_{\frac{p}{2}}}^\frac12 \\
  &  & \kl  \noo xx^*\rrm_{S_{\frac{p}{2}}}^\frac12 \lel
   \noo x\rrm_p \pl .
 \mel
By homogeneity of   $L_{2,r_p}$, this implies $\noo
WL_a\rrm_{cb}\le \noo W\rrm \noo a\rrm_r$.\qed

\begin{cor}\label{ext}
If $T:X_{p,c_p}(\sigma)\to C_p$ is completely bounded,
then $T$ admits a cb-extension to $\ell_p\oplus_p
C_p$.
\end{cor}

\begin{prop} \label{inf}
If $2<p<\8$ and $\N$ is a finite von Neumann
algebra, then $\ell_p(C_p)$ is not cb-isomorphic to
a subspace of $C_p\oplus_p R_p(L_p(\N))$.
\end{prop}

{\bf Proof:} Let $2< r\le \8$ such that
$\frac12=\frac1p+\frac1r$. Let
$T=(T^{(1)},T^{(2)}):\ell_p(C_p)\to C_p\oplus_p
L_p(\N)\oplus_p R_p(L_p(\N))$ be a complete contraction
and $T^{-1}:rg(T)\to \ell_p(C_p)$ be a completely
bounded inverse with $\noo T^{-1}\rrm_{cb} \le C$.
We consider the complete contraction
$T_1:\ell_p(S_p)\to C_p$ defined by
$T_1(x)=T^{(1)}(P(x))$, $P$ the projection onto the
columns space. According to Lemma \ref{leftm}, we can
find $a\in \ell_r(S_r)$  and $W:\ell_2(S_2)\to
\ell_2$ such that $T_1(x)=W(xa)$. Let $\rho=(\noo
a(i)\rrm_r)$ and consider the operator
$D_\rho:\ell_p\to \ell_2$. We define the bounded
map $W':\ell_2(\ell_2)\to \ell_2$ by
$W'((x_i))=W((\rho_i^{-1}x_ia_i))$. In particular,
we can find an $n$ such that
 \[ \kla \summ_{k\ge n} \rho_k^r \mer^{\frac1r} \kl
 \frac{1}{2C} \pl .\]
In the following, we use the  spaces $Y_{n}={\rm
span} \{ \summ_k e_k\ten x_k \p |\p k>n, x_k\in
C_p\}$ and deduce
  \for
  \noo T^{(1)}|_{Y_{n}} \rrm_{cb} &\le&  \kla \summ_{k\ge n_1} \rho_k^r
 \mer^{\frac1r} \pl \noo W':\ell_2(\ell_2)\to
 \ell_2\rrm
 \kl
   \frac{1}{2C}  \pl .
  \mel
If $x\in S_p(Y_n)$, we deduce
 \for
  \frac{1}{C} \noo x\rrm_{S_p(Y_n)} &\le&  \noo (id \ten
  T)(x)\rrm_p
 \kl   \noo (id \ten
 T^{(1)}|_{Y_{n}})(x)\rrm_{C_p}  + \noo (id\ten
 T^{(2)}(x)\rrm_p  \\
 & &\le  \frac{1}{2C} \noo
 x\rrm_{S_p(\ell_p(C_p))}  + \noo (id\ten T^{(1)})(x)\rrm_p
  \pl .
 \mel
Thus
\[ \frac{1}{2C} \noo x\rrm_{S_p(Y_n)} \kl
 \noo (id\ten T^{(1)})(x)\rrm_{S_p(R_p(L_p(\N)))} \kl \noo
  x\rrm_{S_p(\ell_p(C_p))} \pl .\]
In particular $C_p$ is cb-isomorphic to a subspace
of $R_p(L_p(\N))$ which contradicts Lemma
\ref{cpp}.\qed

For the convenience of the reader we quote the following lemma which is used
both in the next proposition and in the next section. The lemma is proved in
\cite{J2} and \cite{JX}.
\begin{lemma}\label{ce}
Let $\M\subset \N$ be von Neumann algebras, $\phi$ a
faithful normal state on $\N$ and $\E:\N \to \M$ a
faithful conditional expectation such that
$\phi|_\M\circ \E=\phi$.  Let $D\in L_1(\M)$ be the
density of $\phi$.
\begin{enumerate}
 \item[i)] If $\frac1r+\frac1s=\frac1p\ge 1$, then
  $\E$ induces a contractive map
 $\E_p:L_p(\N)\to L_p(\M)$ such that
  \[ \E_p(axy) \lel a\E(x)b \]
for all L $a\in L_r(\M)$, $b\in L_s(\M)$ and $x\in
\N$.
\item[ii)] Let $1\le p,p'\le \8$ with $\frac1p+\frac{1}{p'}=1$ and $L_p(\N,\E)$ be
the completion of $\{ aD^{\frac1p}\p | a  \pl
\phi\mbox{ - analytic}  \}$ with respect to the norm
 \[ \noo aD^{\frac1p}\rrm_{L_p(\N,\E)} = \noo
 D^{\frac1p} E(a^*a)D^{\frac1p}
 \rrm_{\frac{p}{2}}^{\frac12} \pl .\]
For $p=\8$, we take the  closure with respect to
strong topology, then
 \[ L_p(\N,\E)^* \lel L_{p'}(\N,\E) \]
and the duality is given by the trace on $\M$.
\item[iii)] Let $1\le p'\le 2\le p\le \8$ with $\frac1p+\frac{1}{p'}=1$, then
 \[ \noo x\rrm_{L_p(\N,\E)} \kl \noo x\rrm_p \]
for all $x\in L_p(\M)$ and
  \[ \noo x\rrm_{p'} \kl \noo x\rrm_{L_{p'}(\N,\E)} \]
 for all $x\in L_{p'}(\N,\E)$.
\end{enumerate}
\end{lemma}

\begin{prop}\label{hcmod}
For every separable subspace $W$ of $\prod_\U
C_p(L_p(0,1))$ there is a commutative von Neumann
algebra $\N$ such that $W$ is completely
isometrically isomorphic to a subspace of
$C_p(L_p(\N))$. If in addition $W$ is
cb-complemented, then $W$ can be assumed
cb-complemented  in  $C_p(L_p(\N))$. The same holds
with $R_p$ replaced by $C_p$.
\end{prop}

{\bf Proof:}  Let us  consider the commutative von Neumann
algebra $\N=(\prod_\U L_1)^*$. Let $\iota:\prod_\U L_1\to
\prod_\U L_1(S_1)$ be the canonical inclusion map, given
coordinatewise by $\iota((f(i))=e_{00}\ten f(i)$ Let
$q_0=(e_{00}\ten 1)$ be the projection onto the first corner.
Obviously $q\le q_0$ and $\E=\iota^*:\prod_\U L_1(S_1)^*\to \N$
defines a conditional expectation. Let $\M=\prod_\U L_1(S_1)^*$
and consider the space
  \[ \noo x\rrm_{S_q^n(L_q(\M,\E))}  \lel \noo (id \ten \E)(x^*x)^{\frac12} \rrm_{S_q(L_q(\N))}  \]
defined on the space of elements $yd^{\frac1q}$,
$d\in L_1(\N)$, $y\in L_q(N)$.    According to
Lemma \ref{ce}, we  have
 \[ L_{p'}(\M,\E)^* \lel L_p(\M,\E) \]
completely isometrically. Obviously, the inclusion
map $T:\prod_\U C_{p'}(L_{p'}(0,1))\to
L_{p'}(\M,\E)$ is completely isometric and
therefore by duality $\prod_\U C_p(L_p(0,1))$ is
completely contractively complemented in
$L_p(\M,\E)$. Given an element $x\in S_p^m(W)$, we
see that
 \[ \noo x\rrm_p^2  \lel \noo x^*x\rrm_{\frac{p}{2}} \lel
 \noo x^*x\rrm_{S_{\frac{p}{2}}^m[L_p(\N)]}  \pl . \]
Since $\bigcup_m S_p^m[W]$ is separable, we can
find a density $D\in L_1(\N)$ such that
 \[ x_{ij}^*x_{ij} \kl C(x) \pl  D^{\frac1p} \pl \]
for all $x=(x_{ij})_{ij=1}^m$ in a countable dense subset $\Delta$
of $\bigcup_m S_p^m[W]$. Multiplying with the support projection $q$
of $D$, we can work in $\N q$. For every coordinate $y=x_{ij}$,
$x=(x_{ij})\in \Delta$, we consider the polar decomposition
 \[ y \lel ub \pl .\]
Using Raynaud's isomorphism \cite{Re}, we see that $b\in
L_p(q\N q)$. Let $\N_1$ be a separable subalgebra
generated by the elements $b=b_{ij}(x)$, $x\in
\Delta$. Let $\M_1$ be a separable subalgebra
containing  by the polar decompositions
$u=u_{ij}(x)$, $x\in \Delta$,  such that there
exists a conditional expectation $\E_1:wcl(\M_1)\to
\N_1$ leaving $\phi$ invariant. Clearly, $W$ is
still a (cb-complemented) subspace of $L_p(\M_1,\E)$
and we can consider the right $\N_1$ module $F$
generated by $M_1$ and $\N_1$. According to
\cite{JX}, $L_p(\M_1,\N_1)$ is completely
contractively complemented in $C_p(L_p(\N_1))$ and
therefore the assertion is proved.\qed

\begin{cor} \label{unifembed}
If $2< p <\8$ and  $F$ is a quotient of
$R_p(L_p(0,1))$, then $\ell_p^n(C_p^n)$ does not embed uniformly
into $ C_p\oplus_p F$.
\end{cor}

{\bf Proof:} Assume to the contrary, we can find
$T_n=(T_n^{(1)},T_n^{(2)}):\ell_p^n(C_p^n)\to
 C_p\oplus_p F$ such that
 \[ \noo T_n\rrm_{cb} \kl 1 \quad \mbox{and} \pl \noo T_n^{-1}\rrm_{cb} \kl C \pl .\]
Let $\U$ be a free ultrafilter on the natural numbers and define
 \[ T: \ell_p(C_p)\to \prodd_\U C_p \oplus_p
 \prodd_\U F \pl ,\]
by $T(x)\lel ((T_n^{(1)}(x))_{\nen},
(T_n^{(2)}(x))_{\nen})$. This is well-defined
because $\bigcup_n \ell_p^n(C_p^n)$ is norm dense
in $\ell_p(C_p)$. Moreover, for $x\in
S_p^m(\ell_p^n(C_p^n))$, we have
 \begin{equation*}
  \noo (id \ten T) (x)\rrm = \lim_{n'>n} \noo id \ten T_{n'}(x)\rrm_{S_p(\ell_p)\oplus_p S_p(C_p)}
  \sim_C \noo x\rrm_{S_p^m((\ell_p^n(C_p^n))} \pl .
\end{equation*}
Let us denote the first component by $T^{(1)}$ and
the second by $T^{(2)}$.  We note that $\prod_\U F$
is a quotient space of $\prod_\U R_p(L_p(0,1))$.
Denote the quotient map by $q$. Then we can find a
separable subspace $Y\subset \prod_\U
R_p(L_p(0,1))$ such that the image of $T^{(2)}$ is
cb isomorphic to $q(Y)$. According to Proposition
\ref{hcmod}, we can assume that $Y$ is contained in
$R_p(L_p(\N))$ for some commutative von Neumann
algebra $\N$. Moreover, $\prod_\U C_p$ is a
homogeneous   Hilbert space which carries  the
$C_p$ structure. Thus every separable subspace is
completely isometric to $C_p$. Therefore, we can
find an embedding of $\ell_p(C_p)$ in $C_p\oplus_p
Y/ker(q)$. Following the argument in Proposition
\ref{inf}, we see that for the first component
$T^{(1)}$ and every $\eps>0$ there exists an $n$
such that $\noo T^{(1)}|_{\{(x_k) \p |
x_1=x_1=\cdots=x_{n}=0\}}\rrm_{cb} \kl \eps$. Thus
$C_p$ will be $cb$-isomorphic to a subspace of a
quotient  of $R_p(L_p(0,1))$. This contradicts Lemma
\ref{dist}. \qed

\begin{theorem} \label{exclude}
Let $\si$ tend to $0$ and such that
for all $\nen$ there are subset $A_n$ of cardinality $n$ such that
$\si_i=\al_n$ for $i\in A_n$   and
 \[  \lim_n n^{\frac1r}\al_n \lel \8 \pl .\]
Then $X_{p,c_p}(\si)$ does not admit a cb factorization
through $C_p\oplus_p F$, $F$ a quotient of a
subspace of $\prod_\U R_p(L_p(0,1))$.
\end{theorem}

\hz

{\bf Proof:}  Assuming in the contrary we can write
$id=T+S$, where $T$ factors through a quotient $F$ of $\prod_\U
R_pL_p(0,1)$ and $S$ factors through $C_p$. We
denote by $Q$ the projection onto the $C_p$
coordinate in $X_{p,c_p}(\si)\subset \ell_p\oplus_p
C_p$. Using Lemma \ref{ext},  we can decompose
$S=S_1+S_2$, such that $S_1:\ell_p\to X_{p,c_p}$ is
a completely bounded operator and $S_2: C_p\to
X_{p,c_p}$ is completely bounded.  For a fixed index
$i\in I$ we consider
 \for
 (e_i,\si_ie_i) &=& S(e_i,\si_i e_i) +
 T(e_i,\si_ie_i)\lel \si_i S_2(0,e_i) + S_1(e_i,0)
 +T(e_i,\si_ie_i)\pl .
 \mel
Thus
 \[ 1\le \noo   S_1(e_i,0) +T(e_i,\si_ie_i)\rrm + \si_i \noo S_2\rrm \pl .\]
Hence for $i\ge i_0$ we get $\si_i \noo S_2\rrm \le \frac12$ and
therefore
 \[ \frac12 \kl  \noo   S_1(e_i,0) +T(e_i,\si_ie_i)\rrm  \pl .\]
Let us write
  \[ S_1(e_i,0) +T(e_i,\si_ie_i)\lel (y,\si y) \pl .\]
We have the following alternative:   If $\noo y\rrm_p \le  \noo
y\si\rrm_2$, then
 \for
 \frac12 \le \kla \noo y\rrm_p^p +\noo y\si\rrm_2^p\mer^{\frac1p}
 \kl 2\noo y\si\rrm_2 \pl .
 \mel
Hence
  \[ \frac14 \kl \noo y\si\rrm_2 \pl .\]
If $\noo y\si\rrm_2\le \noo y\rrm_p$, we get
 \[ \frac14 \kl \noo y\rrm_p \]
and thus
 \[ \frac{\si_i}{4}  \noo y\rrm_p  \kl \noo y\si\rrm_p \kl \noo y\si\rrm_2 \pl .\]
In both cases we deduce
 \[ \frac{\si_i}{4} \kl \noo QS_1(e_i,0)+   QT(e_i,\si_ie_i)\rrm_2 \pl .\]
Now we decompose $QT=T_1+T_2$, $T_1$ acting on
$\ell_p$ and $T_2$ acting on $C_p$    according to
Lemma \ref{ext}. Let $\nen$ to be determined later
and let us assume that $\si_i=\al_n$ is constant on
a set $A_n$ of cardinality $n$. Let us recall that
 \[ \kla \summ_i \noo QS_1(e_i)\rrm_2^r  \mer^{\frac1r} \kl \noo
 QS_1\rrm \kl C_1 \]
and
  \[ \kla \summ_i \noo T_1(e_i)\rrm_2^r  \mer^{\frac1r} \kl \noo
 T_1\rrm \kl C_2 \pl . \]
Thus we get for $C_3=\noo T_2\rrm $
 \for
 & &\frac{\al_n n^{\frac1r} }{4} \le \kla  \summ_{i\in A_n} \noo QS_1(e_i,0)+
 QT(e_i,\si_ie_i)\rrm_2^r \mer^{\frac1r} \kl
 C_1 + C_2 + \kla \summ_{i\in A_n} \noo T_2(0,\si_i e_i)
 \rrm_2^{r} \mer^{\frac1r}  \\
& & \kl C_1 + C_2 + \kla \summ_{i\in A_n, \noo T_2(0,e_i)\rrm \kl
\frac{1}{16}}   \!\!\!   \noo T_2(0,\si_i e_i)
 \rrm_2^{r} \mer^{\frac1r}
  +
 \kla \summ_{i\in A_n, \noo T_2(0,e_i)\rrm > \frac{1}{16}} \!\!\!   \noo T_2(0,\si_i e_i)
 \rrm_2^{r} \mer^{\frac1r} \\
 & & \kl  C_1+C_2  + \al_n  \frac{1}{16}  n^{\frac1r} +
 \al_n C_3   {\rm card}\{ i\in A_n \p | \p \noo T_2(0, e_i)  \rrm >\frac{1}{16}
 \} \pl .
 \mel
Hence for $n$ so large that $8(C_1+C_2)\le \al_n
n^{\frac1r}$ we get
 \[ \frac{1}{16 C_3} n^{\frac1r}  \kl
  {\rm card}\{ i\in A_n \p | \p \noo T_2(0, e_i)  \rrm
  >\frac{1}{16}\} \pl .\]
Hence we can find a subset $B_n$ of cardinality
$\frac{n}{C_3^r16^r}$ such that for all $i\in B_n$ we have
 \[ \noo T_2(0, e_i)  \rrm_2  >\frac{1}{16} \pl .\]
Now we consider the map $w:\ell_2(B_n)\to \ell_2$
defined by $w(e_i)=T_2(0, e_i)$. Defining
$\delta=C_3^{-1}32^{-2}$  and $n'=card{B_n}$
we deduce for the approximation numbers of $w$
 \for
 \frac{1}{16} \pl \sqrt{n'} &\le& \pi_2(w) \kl
 \kla \summ_{k=1}^{n'} a_k(w)^2 \mer^{\frac12}
  \kl  \sqrt{\delta}  \sqrt{n'}\noo T_2\rrm  + a_{\delta n'}(w)
  \sqrt{n'}\\
 &\le& \frac{1}{32} \sqrt{n'} + a_{\delta n'}(w) \sqrt{n'} \pl .
 \mel
Therefore with $\delta'=C_3^{-r}16^{-r}$ we obtain
that
 \[ \frac{1}{32} \kl a_{\delta n'}(w) \lel a_{\delta \delta' n}(w)  \pl
 .\]
Let $u:\ell_2(B_n)\to C_p \cong \ell_2$ be defined
by $u(e_i)=QT(e_i,\si_ie_i)$. In order to obtain a
lower estimate for a proportional approximation
number of $u$ we observe
 \[ \al_n w(e_i) \lel T_2(0,\si_ie_i) \lel QT(e_i,\si_ie_i) -
 T_1(e_i,0)  \lel u(e_i)-T_1(e_i,0)  \pl .\]
Since $T_1$ is bounded on $\ell_p$, the  map
$T'_1:\ell_2\to \ell_2$  defined by   $e_i\mapsto
T_1(e_i,0)$ factors through the inclusions map
$id_{2,p}:\ell_2\to \ell_p$
 \[ \al_n w-u \lel T_1 id_{p,2} \pl ,\]
Let us recall   a result of Carl on the Weyl numbers
of $id_{p',2}$
 \[ k^{\frac{1}{r}} x_k(id:\ell_{p'}\to \ell_2) \kl c_0 \pl .\]
Therefore we have
 \for
 \frac{\al_n}{32} &\le& a_{\delta \delta' n}(\al_n w) \lel
 a_{\delta \delta' n}(u+\al_nw-u) \\
 &\le & a_{\frac{\delta \delta'}{2} n}(u) +
 a_{\frac{\delta \delta'}{2} n}(T_1id_{p,2})
  \lel  a_{\frac{\delta \delta'}{2} n}(u) + \kla \frac{2n}{\delta
 \delta'} \mer^{-\frac1r} c_0 \noo T_1\rrm \pl .
 \mel
Hence for $n$ large enough such that $n^{\frac1r}\al_n
\gl\frac{128 c_0 \noo T_1\rrm}{\delta \delta'}$ we obtain
 \[ \frac{\al_n}{64}   \kl a_{\frac{\delta \delta'}{2} n}(u)
 \pl .\]
It follows that we can find an linear map $W:\ell_2\to
\ell_2$ and a $k=\frac{\delta \delta'}{2} n$
dimensional subspace $H\subset \ell_2(B_n)$ such
that $\noo W\rrm\le 64 \al_n^{-1}$ and $WQTP_H \lel id_H $.

Note that cb norm of the identity mapping
$id:C_p\to X_{p,c_p}$ is completely contractive and
thus we obtain
 \[ id_H \lel WQTidP_H  \pl .\]
According to our assumption  $T=w_1v_1$ where
$v_1:X_{p,c_p}(\si) \to F$, $w_1:F\to
X_{p,c_p}(\si)$ and $F$ is a quotient to a
subspace of $\prod_\U R_p(L_p(0,1))$. We deduce from
Lemma \ref{dist2} that
 \for
  \frac{\delta \delta'}{2}n^{\frac1r} &=&
  k^{\frac1r} \kl c_p \pl  \inf_{E\in QS(\prod_\U R_p(L_p(0,1)))} d_{cb}(C_p^k,E)  \\
  &\le&  c_p \pl \noo W\rrm_{cb}  \noo v_1\rrm_{cb} \noo
  w_1\rrm_{cb}
  \kl   \al_n^{-1} c_p \pl \noo v_1\rrm_{cb} \noo
  w_1\rrm_{cb} \pl .
 \mel
Using once more $\lim_n n^{\frac1r} \al_n =\8$, we
get a contradiction and the assertion is
proved.\qed

\begin{theorem}
\label{thm2.20n}
If $V \subseteq \ell_p \oplus_p C_p \oplus_p R_p$ is a rectangular $\cO
\cL_p$-space, then there exists an increasing sequence $(X_j)$ of finite
dimensional subspaces of $V$ with dense union and non-negative integers
$k_j$, $m_j$, $n_j$ and a constant K so that

\begin{equation}
\label{eq2.1ny}
d_{cb}(X_j,\ell_p^{k_j} \oplus_p C_p^{n_j} \oplus_p R_p^{m_j}) \le K \quad
\mbox{for all $j\in \bN$}.
\end{equation}

In particular $V$ is cb-isomorphic to a cb-complemented subspace of
$L_p(0,1) \oplus_p C_p \oplus_p R_p$.

If $V \subseteq \ell_p \oplus_p C_p$, the $R_p$-terms in (\ref{eq2.1ny})
disappear and $V$ is cb-isomorphic to a cb-complemented subspace of
$L_p(0,1) \oplus_p C_p$. Similarly if $V \subseteq \ell_p \oplus_p R_p$.
\end{theorem}

\bproof
Since $V$ is a rectangular $\cO \cL_p$-space there is an increasing
sequence $(X_j)$ of finite dimensional subspace with dense union and
number $k(j)$, $n_j(i)$ and $m_j(i)$  and a constant $K_1$ so that
\begin{equation*}
d_{cb}(X_j,(\oplus_{i=1}^{k(j)} S_p^{n_j(i),m_j(i)})_p) \le K_1 \quad
\mbox{for all $j\in \bN$}.
\end{equation*}
For every $n\in \bN$ we define
\begin{equation*}
h(n) = \sup\{m_j(i) \mid n_j(i) \ge n \}.
\end{equation*}
If $h(n) \ge n$ for all $n\in \bN$, clearly $(S_p^n)$ embeds cb-uniformly
into $V$ and hence $S_p$ is isomorphic to a subspace of an ultrapower of
$\ell_p \oplus_p C_p \oplus_p R_p$ which is a Banach lattice of cotype
p. This contradicts \cite[Theorem 2.1]{P2}. Hence there is an $n_0\in
\bN$ so that $h(n_0) < n_0$. If $n_j(i) \le n_0$, then
\begin{equation*}
d_{cb}(S_p^{n_j(i),m_j(i)}, \ell_p^{n_j(i)}(R_p^{m_j(i)}) \le
n_0^{\frac{1}{r}}
\end{equation*}
and if $n_j(i)\ge n_0$, then $ m_j(i) < n_0$ and hence
\begin{equation*}
d_{cb}(S_p^{n_j(i),m_j(i)}, \ell_p^{m_j(i)}(C_P^{n_j(i)})) \le n_0^{\frac{1}{r}}.
\end{equation*}

We can therefore find a constant $K_2$ and numbers $k_j^{'}$,
$n_j^{'}(i)$ and $m_j^{'}(i)$ so that
\begin{equation*}
d_{cb}(X_j, (\oplus_{i=1}^{k_j^{'}} C_p^{n_j^{'}(i)})_p \oplus_p
(\oplus_{i=1}^{k_j^{'}} R_p^{m_j^{'}(i)})_p) \le K_2 \quad \mbox{for all
$j\in \bN$}.
\end{equation*}

For every $n$ and $j$ we put $A_j(n) = \{i \le k_j^{'} \mid n_j^{'}(i) \ge
n \}$ and $f(n) = \sup_j |A_j(n)|$. If $f(n) \ge n$ for all $n \in \bN$,
then clearly $(\ell_p^n(C_p^n))$ embeds cb-uniformly into $V \subseteq
\ell_p \oplus_p C_p \oplus_p R_p$ which contradicts Corollary
\ref{unifembed}. Hence there is an $n_0$ so that $|A_j(n_0)|<n_0$ for
all $j\in \bN$. For every $j$ we then get
\begin{eqnarray*}
d_{cb}((\oplus_{i \in A_j(n_0)}C_p^{n_j^{'}})_p,C_p^{\sum_{i \in A_j(n_0)}
n_j^{'}(i)}) & \le & n_0^{\frac{1}{r}}\\
d_{cb}(\oplus_{i \notin A_j(n_0)} C_P^{n_j^{'}(i)}, \ell_p^{\sum_{i \notin
A_j(n_0)} n_j^{'}(i)}) & \le & n_0^{\frac{1}{r}}.
\end{eqnarray*}

Treating the $R_p$-terms in the same way we obtain that there is a
constant $K$ and numbers $k_j$, $n_j$ and $m_j$ so that
\begin{equation*}
d_{cb}(X_j, \ell_p^{k_j} \oplus_p C_p^{n_j} \oplus_p R_p^{m_j}) \le K \quad
\mbox{for all $j \in \bN$}
\end{equation*}
which proves formula (\ref{eq2.1ny}). Using an ultraproduct construction as in \cite[Section 10.3]{ER2} we
deduce that there is an ultrafilter $\U$ so that $V$ is cb-complemented
in $\prod_{\U}\ell_p \oplus_p \prod_{\U} C_p \oplus_p \prod_{\U}
R_p$. Since $\prod_{\U} \ell_p$ is cb-isometrically isomorphic to
$L_p(\cN)$ for some commutative $\cN$ and $C_p$ and $R_p$ are
homogeneous, the separability of $V$ gives that it is cb-complemented in
$L_p(\cN_1) \oplus_p C_p \oplus_p R_p$ with $(\cN_1)_*$
separable. Decomposing $\cN_1$ into discrete and continuous parts we
get that $L_p(\cN_1)$ is cb-contractively complemented in $L_p(0,1)$ and
hence $V$ is cb-isomorphic to a cb-complemented subspace of $L_p(0,1)
\oplus_p C_p \oplus_p R_p$.

Since $(R_p^n)$ does not embed cb-uniformly into $\ell_p \oplus_p C_p$ by
Lemma \ref{dist}, it is readily seen that if $V \subseteq \ell_p \oplus_p
C_p$, then the $R_p$-components disappear in the argument above and the
ultraproduct construction gives that $V$ is cb-isomorphic to a
cb-complented subspace of $L_p(0,1) \oplus_p C_p$.
\eproof

As a corollary we obtain

\begin{theorem}
\label{prop2.20}
If $\sigma$ satisfies (\ref{eq1.1}) and (\ref{eq1.2}), then
the spaces $X_{p,c_p}(\sigma)$ and $X_{p,r_p}(\sigma)$are not rectangular
$\cO\cL_p$ spaces.
\end{theorem}

\bproof
Assume that $X_{p,c_p}(\sigma)$ is a rectangular $\cO\cL_p$-space. Theorem \ref{thm2.20n}
then gives that it is cb-complemented in $L_p(0,1)\oplus_p C_p$. By Theorem \ref{cor1.1a} we can without loss of generalty assume
that $\sigma$ satisfies the additional assumptions in
Theorem \ref{exclude} and hence this theorem
yields a contradiction.\qed

\begin{theorem} \label{rolp}
If $\sigma$ satisfies (\ref{eq1.1}) and (\ref{eq1.2}) and
 \[ V\in\{ R_p\oplus_p X_{p,c_p}(\si), \ell_p(R_p)\oplus_p
 X_{p,c_p}(\si), X_{p,r_p}(\si) \oplus_p X_{p,c_p}(\si)\}\pl ,\]
then $V$ is not a rectangular  $\cO\cL_p$ space.
\end{theorem}

{\bf Proof:} Let us assume $V= \ell_p(R_p)\oplus_p X_{p,c_p}(\sigma)$. The proof of Theorem \ref{thm2.20n} shows that $V$ is cb-complemented in $C_p \oplus_p \prod_{\U}R_p(\ell_p)$ which contradicts Theorem \ref{exclude} since $X_{p,c_p}(\sigma)$ is cb-complented in $V$. The other cases follow directly from Theorem \ref{thm2.20n}.
\eproof

\begin{prop}
\label{propultra}
Assume that $\sigma$ satisfies (\ref{eq1.1}) and (\ref{eq1.2})and let $\U$ a free ultrafilter on the integers.
\begin{enumerate}
\item[(i)] If $V\in \{X_{p,c_p}(\si),R_p\oplus_p X_{p,c_p}(\si),
X_{p,r_p}(\si) \oplus_p X_{p,c_p}(\si)\}$, then $\ell_p(R_p)\oplus_p
X_{p,c_p}(\si)$ does not embed into $\prod_\U V$.
\item[(ii)] $X_{p,r_p}(\si) \oplus_p X_{p,c_p}(\si)$ is not cb-isomorphic to a cb-complemented subspace of $\prod_\U (R_p\oplus_p X_{p,c_p}(\si) ) $.
\end{enumerate}

In particular the spaces $\{X_{p,c_p}(\si), R_p\oplus_p
X_{p,c_p}(\si), X_{p,r_p}(\si) \oplus_p X_{p,c_p}(\si),
\ell_p(R_p)\oplus_p X_{p,c_p}(\si)\}$ are mutually not
cb-isomorphic.
\end{prop}

{\bf Proof:} To prove the assertion $(i)$, we observe that
$V\subset \ell_p \oplus_p C_p \oplus_p R_p$. Thus the assertion
follows from the row version of Corollary \ref{unifembed}. In order to get
$(ii)$ we note that $R_p\oplus_p X_{p,c_p}(\sigma)$ is complemented
in $R_p\oplus_p L_p([0,1];C_p)$. According to Proposition
\ref{hcmod} a separable complemented subspace of $\prod_\U
R_p\oplus_p L_p([0,1];C_p)$ is cb-complemented  in $R_p\oplus_p
C_p(L_p(\N))$ for a commutative $\N$. But the row version of
Theorem \ref{exclude} excludes this for $X_{p,r_p}(\si)$.\qed

\begin{rem}
\label{remultra}
If $W \in \{\ell_p(R_p),\ell_p(R_p)\oplus_p X_{p,c_p}(\si), R_p\oplus_p X_{p,c_p}(\si),
X_{p,r_p}(\si) \oplus_p X_{p,c_p}(\si)\}$, then $W$ contains $R_p$
cb-somorphically which does not cb-embed into an ultrapower of
$L_p([0,1];C_p)$. However, $X_{p,c_p}(\sigma) \subseteq L_p([0,1];C_p)$
and hence $W$ does not cb-embed into an ultrapower of
$X_{p,c_p}(\sigma)$.

Consequently none of the spaces above nor those from Proposition \ref{propultra}
can be paved with local pieces of any of the others
except for trivial reasons. It is easily seen that
we can also add  $\ell_p(C_p)\oplus_p X_{p,r_p}(\sigma)$ and
the rectangular $\cO\cL_p$ space
$\ell_p(C_p)\oplus_p X_{p,c_p}(\sigma)$  to this list.
\end{rem}

At the end of this section we want to compare the space
$X_{p,c_p}(\si)$, $X_{p,r_p}(\rho)$  with their intersection in
interpolation sense. Let $2<p<\infty$ and let $\sigma = (\sigma_n)$ and $\rho = (\rho_n)$ be two positive sequences. In
analogy with the spaces defined in chapter 1 we let the space
$X(\sigma,\rho)$ be the subspace of $S_p \oplus_p C_p \oplus_p R_p$
defined as the closed linear span of the sequence $\{e_{nn} \oplus_p
\sigma_n e_{n1} \oplus_p \rho_n e_{1n}\}$. Note that $X(\sigma,\rho)$ is the
interpolation space $X_{p,c_p}(\sigma) \cap X_{p,r_p}(\rho)$. We shall
show that if $\sigma$ and $\rho$ satisfy (\ref{eq1.1}) and (\ref{eq1.2}),
then $X_p(\sigma,\rho)$ is a rectangular $\cO\cL_p$-space if and only if
it is cb-isomorphic to $X_p(\sigma)$, $X_p(\sigma)\oplus_p C_p$, $X_p(\sigma) \oplus_p R_p$ or $X_p(\sigma) \oplus_p C_p \oplus_p R_p$. We first investigate the
space $X_p(\alpha,\beta)$ where $\alpha > 0$ and
$\beta>0$ are constants. We have:

\begin{proposition}
\label{prop2.24n}
There is a constant $K=K(p)$ so that if $T$ is a cb-isomorphism of
$X_p(\alpha,\beta)$ into $L_p(0,1)\oplus_p C_p \oplus_p R_p$ and
$P$ is a cb-projection of $L_p(0,1) \oplus_p C_p \oplus_p R_p$ onto
$T(X_p(\alpha,\beta))$, then either
\begin{equation}
\label{eq2.1n}
\max(\alpha,\beta) \le K\|T\|_{cb}\|T^{-1}\|_{cb}\min(\alpha,\beta)
\end{equation}
or
\begin{equation}
\label{eq2.2n}
\frac{1}{2\min(\alpha, \beta)} \le K \|P\|_{cb}\|T\|_{cb}\|T^{-1}\|_{cb}
\end{equation}

If $T$ is a cb-isomorphism of $X_{p,c_p}(\alpha)$ into $L_p(0,1) \oplus_p C_p$ and $P$ is a cb-projection of $L_p(0,1) \oplus_p C_p$ onto $T(X_{p,c_p}(\alpha))$, then
\begin{equation}
\label{eq2.2ny}
\frac{1}{2 \alpha} \le K \|P\|_{cb}\|T\|_{cb}\|T^{-1}\|_{cb}.
\end{equation}
Similarly for $X_{p,r_p}(\sigma)$.
\end{proposition}

\bproof Let us assume that $\beta \le \alpha$ (the other case can be
proved similarly), let $Q_1$ be the natural projection of
$L_p(0,1)\oplus_p C_p\oplus_p R_p$ onto $L_p(0,1)$ and $Q_2$ the natural
projection of $L_p(0,1)\oplus_p C_p \oplus_p R_p$ onto $C_p\oplus_p R_p$.
If $(f_n)$ denotes the canonical basis of $X_p(\alpha,\beta)$, we
put $h_n = Q_1Tf_n$ for all $n\in \bN$. Since $f_n \to 0$ weakly, so
does $(h_n)$ and we can therefore extract a martingale subsequence
of $(h_n)$ and then use the argument in \cite{JXII} to extract a
further subsequence, still called $(h_n)$, so that there exist
constants $K_1 = K_1(p)\ge 1$, $b_1\ge 0$ and $b_2\ge 0$ so that
\begin{equation*}
\|\sum_k a_k h_k\|_{S_p[L_p(0,1)]}
\sim_{K_1}
\max\{b_1\sum_k\|a_k\|_{p}^p, b_2\|
(\sum_k a_k^*a_k)^{\frac{1}{2}}\|_{S_p},
b_2 \| (\sum_k a_ka_k^*)^{\frac{1}{2}}\|_{S_p}\}
\end{equation*}
for all finite sequences $(a_k) \subseteq S_p$. Plugging in the vectors
$a_k = e_{1k}$ we get for every $n\in \bN$ that
\begin{equation*}
\max(b_1 n^{\frac{1}{p}}, b_2 n^{\frac{1}{p}}, b_2 n^{\frac{1}{2}}) \le
K_1 \|T\|_{cb}\max(n^{\frac{1}{p}}, \alpha n^{\frac{1}{p}}, \beta
n^{\frac{1}{2}})
\end{equation*}
which implies that $b_2 \le K_1\beta$.

As in Corollay \ref{ext} there is a constant $K_2$ only depending on p so
that the operator $Q_2T$ has a cb-extension $S: S_p \oplus_p C_p \oplus_p
R_p \to C_p \oplus_p R_p$ with $\|S\|_{cb} \le K_2\|T\|_{cb}$. Hence we have
for all $n\in \bN$:
\begin{equation*}
Tf_n = h_n + Se_{nn} + \alpha Se_{n1} + \beta Se_{1n}.
\end{equation*}
By \cite{R} $\sum_{n=1}^{\infty} \|Se_{nn}\|^r < \infty$ and if $Q$ denotes the canonical projection of $S_p \oplus_p C_p \oplus_p R_p$ onto $R_p$ we find
by that the operator $QT^{-1}PS|C_p$ is (r,2)-summing and therefore also\\
$\sum_{n=1}^{\infty}\|QT^{-1}PSe_{n1}\|^r < \infty$. In particular we can
find a $n_0 \in \bN$ so that:
\begin{equation}
\label{eq2.3n}
\|T^{-1}PSe_{nn}\| + \frac{\alpha}{\beta} \|QT^{-1}PSe_{n1}\| \le \frac{1}{4}
\end{equation}
for all $n \ge n_0$. If $(F_n)$ denotes the biorthogonal system to
$(f_n)$, then clearly $|F_n(T^{-1}PSe_{n1})| \le
\frac{1}{\beta}\|QT^{-1}PSe_{n1}\|$ and hence (\ref{eq2.3n}) gives
that
\begin{eqnarray*}
1 &\le&  |F_n(T^{-1}Ph_n)| +
\beta |F_n(T^{-1}PSe_{1n}| + \frac{1}{4}  \\
 &\le& |F_n(T^{-1}Ph_n)| + K_2\beta \|P\|_{cb}\|\|T\|_{cb}\|T^{-1}\|_{cb}+\frac{1}{4}.
\end{eqnarray*}
for all $n\ge n_0$. If we now assume that $K_2
\beta\|T\|_{cb}\|T^{-1}\|_{cb}\|P\|_{cb} < \frac{1}{2}$, then by the
above $\frac{1}{4}\le |F_n(T^{-1}Ph_n)|$ for all $n\ge n_0$.

By interpolation there exists a constant $K_3 = K_3(p)$ so that if $U$
denotes the diagonal of $T^{-1}P|[h_n]$ with respect to the bases $(f_n)$ and
$(h_n)$, then $U$ is cb-bounded with $\|U\|_{cb} \le
K_3\|T^{-1}\|_{cb}\|P\|_{cb}$ and hence for all $(a_k) \subseteq S_p$
and all $n\ge n_0$ we get:
\begin{eqnarray*}
\frac{1}{4} \|\sum_{k=n_0}^n a_k \otimes f_k\| & \le &
\|U(\sum_{k=n_0}^n a_k \otimes h_k)\|_{S_p[L_p(0,1)]} \\
& \le & K_3 \|T^{-1}\|_{cb}P\|_{cb}\|\sum_{k=n_0}^n a_k \otimes
h_k\|_{S_p[L_p(0,1)]}.
\end{eqnarray*}

If we plug in the vectors $a_k = e_{k1}$ in this inequality we get that
\begin{eqnarray*}
& & \frac{1}{4} \max\{(n-n_0)^{\frac{1}{p}}, \alpha
(n-n_0)^{\frac{1}{2}},
\beta (n-n_0)^{\frac{1}{p}}\}  \\
& & \kl
K_1K_3\|T\|_{cb}\|T^{-1}\|_{cb}\|P\|_{cb}\max\{b_1(n-n_0)^{\frac{1}{p}},
b_2 (n-n_0)^{\frac{1}{2}}, b_2 (n-n_0)^{\frac{1}{p}}\}
\end{eqnarray*}
and therefore $\alpha \le K_1K_3\|T\|_{cb}\|T^{-1}\|_{cb}\|P\|_{cb}b_2
\le K_1^2K_3\|T\|T^{-1}\|_{cb}\|P\|_{cb}\beta$. Hence we have proved
the proposition with $K=\max(K_1^2K_3,K_2)$.

To prove the statement for $X_{p,c_p}(\alpha)$ we go through the
argument above, but we omit the $R_p$-coordinate, and  adjust the
sequence $(h_n)$ to the new situation. Then we drop the argument
with the projection $Q$. The first part will then show that $b_2 \le
K_1 \alpha$. If $K_2 \alpha \|T\|_{cb}\|T^{-1}\|_{cb}\|P\|_{cb} <
\frac{1}{2}$, then the second part will show that $\alpha \le K_1K_3
\|T\|_{cb}\|T^{-1}\|_{cb}\|P\|_{cb}b_2$. Hence $(f_n)$ is
cb-equivalent to $(h_n)$ which is a contradiction because
$X_{p,c_p}(\alpha)$ is cb-isomorphic to $C_p$ which does not embed
into $L_p(0,1)$ by Lemma \ref{cpp}. \eproof

We need the following two lemmas:

\begin{lemma}
\label{2.25n}

Let $2\le p < \infty$ and let $\sigma$ and $\rho$ be two sequences so
that there exists a $\delta>0$ and an $\varepsilon > 0$ with $\sigma_n \le
\delta \rho_n$ for all $n \in \bN$ and $\sum_{\sigma_n \le
\varepsilon}\sigma_n^r <\infty$.

If $X_p(\sigma,\rho)$ is cb-isomorphic to a cb-complemented subspace of
$L_p(0,1)\oplus_p C_p \oplus_p R_p$, then
there exist $0\le K,M,N \le \infty$ so that $X_p(\sigma,\rho)$
is cb-isomorphic to $l_p^N \oplus_p (C_p\cap R_p)^M \oplus_p R_p^K$.

If $\rho_n \to 0$, the last two summands do not occur in the above.
\end{lemma}

\bproof Assume that $X_p(\sigma,\rho)$ is a $\cO\cL_p$-space, put
\begin{eqnarray*}
A & = & \{n \in \bN \mid \sigma_n \le \varepsilon\}\\
B & = & \{n \in \bN \mid \sigma_n > \varepsilon \}
\end{eqnarray*}
and let $\sigma_A = \{\sigma_n \mid n \in A\}$ and $\sigma_B = \{\sigma_n
\mid n\in B\}$. In a similar
manner we define $\rho_A$ and $\rho_B$. Clearly we can write
\begin{equation*}
X_p(\sigma,\rho) = X_p(\sigma_A,\rho_A) \oplus X_p(\sigma_B, \rho_B).
\end{equation*}
If $\liminf\rho_A(n) >0$, $X_p(\sigma_A, \rho_A)$ is
cb-isomorphic to $R_p^{|A|}$ (which is cb-isomorphic to $\ell_p^{|A|}$ in
case $A$ is finite). Assume next that $\liminf\rho_A(n) =0$. If
$\rho_A$ satisfies (\ref{eq1.2}), $X_p(\sigma_A,\rho_A)$ is
cb-isomorphic to $X_{p,r_p}(\rho_A)$ which contradicts Theorem
\ref{prop2.20} and hence there is an $\varepsilon_1 >0$ so that
$\sum_{\rho_A(n) \le \varepsilon_1} \rho_A(n)^r < \infty$. We may
without loss of generality assume that $\varepsilon_1 = \varepsilon$ and
can conclude that $X_p(\sigma,\rho)$ is cb-isomorphic to $\ell_p^{|A|}$.
If $n\in B$, $\varepsilon < \sigma_n \le \delta \rho_n$ so that
$X_p(\sigma_B,\rho_B)$ is cb-isomorphic to $(C_p \cap R_p)^{|B|}$.

Summing up we have found that there exist $0\le K,M,N \le \infty$ so that
$X_p(\sigma,\rho)$ is cb-isomorphic to $\ell_p^N \oplus_p (C_p \cap
R_p)^M \oplus_p R_p^K$.
\eproof

\begin{lemma}
\label{2.26n}
Let $2<p< \infty$ and let $\sigma$ and $\rho$ be two sequences so that
$X_p(\sigma,\rho)$ is cb-complemented in $L_p(0,1)\oplus_p C_p \oplus_p R_p$.
Then $\{\rho_n \mid
\sigma_n\ge \varepsilon\}$ does not satisfy (\ref{eq1.2}) for any
$\varepsilon >0$.  The same holds with $\sigma$ and $\rho$
interchanged.
\end{lemma}

\bproof Assume that there is an $\varepsilon > 0$ so that $\{\rho_n
\mid \sigma\ge \varepsilon\}$ satisfies (\ref{eq1.2}).Then it also
satisfies (\ref{eq1.1}) and if $\beta >0$ is arbitrary, we can find
a sequence $(B_k)$ consisting of mutually disjoint finite subsets of
$\bN$ so that
\begin{equation*}
\beta \le (\sum_{n\in B_k, \sigma_n\ge \varepsilon}\rho_n^r)^{\frac{1}{r}} \le 2 \beta.
\end{equation*}
For every $k\in \bN$ we put $\alpha_k = (\sum_{n\in B_k, \sigma_n
\ge \varepsilon}\sigma_n^r)^{\frac{1}{r}}$ and arguing like in
Proposition \ref{thm1.1a} we get that $X_p((\alpha_k),\beta)$ is
cb-complemented in $X_p(\sigma, \rho)$. Clearly $\alpha =
\liminf\alpha_k \ge \varepsilon$ and if we choose a subsequence
$(\alpha_{k_m})$ tending sufficiently fast to $\alpha$ we conclude
that $X_p(\alpha, \beta)$ is cb-complemented in $X_p(\sigma, \rho)$
and hence also in $L_p(0,1)\oplus_p C_p \oplus_p R_p$. This violates
(\ref{eq2.1n}) and (\ref{eq2.2n}) for $\beta$ small enough. \eproof

We are now able to prove:
\begin{theorem}
\label{2.27n}
Let $\sigma$ and $\rho$ be two sequences satisfying (\ref{eq1.1}) and
(\ref{eq1.2}). If $X_p(\sigma, \rho)$ is a rectangular $\cO
\cL_p$-space, then it is cb-isomorphic to $X_p(\sigma)$,
$X_p(\sigma)\oplus_p R_p$, $X_p(\sigma) \oplus_p C_p$ or $X_p(\sigma) \oplus_p
C_p\oplus_p R_p$.

If in addition $\sigma_n \to 0$ and $\rho_n \to 0$, $X_p(\sigma,\rho)$ is
cb-isomorphic to $X_p(\sigma)$.
\end{theorem}

\bproof  If $X_p(\sigma,\rho)$ be a rectangular $\cO\cL_p$-space,
Theorem \ref{prop2.20} shows that it is cb-isomorphic to a
cb-complemented subspace of $L_p(0,1)\oplus_p C_p \oplus_p R_p$. Assume
that for all $\varepsilon >0$ and all $\delta >0$ we have that
$\sum_{\{\rho_n \le \delta \sigma_n, \rho_n \le
\varepsilon\}}\rho_n^r = \infty$. We shall show that this leads to a
contradiction. Let $\delta > 0$ be given arbitrarily, put $A =\{n
\in \bN \mid \rho_n \le \delta \sigma_n \}$ and define $\sigma_A$
and $\rho_A$ as in Lemma \ref{2.25n}. Clearly $\rho_A$ satisfies
(\ref{eq1.1}) and (\ref{eq1.2}). If for some $\varepsilon > 0$
$\sum_{\{\sigma_A(n)\le \varepsilon \}}\sigma_A(n)^r < \infty$, then
also $\sum_{\{\sigma_A(n) \le \varepsilon\}} \rho_A(n)^r < \infty$
and therefore $\{\rho_A(n) \mid \sigma_A(n) > \varepsilon\}$
satisfies (\ref{eq1.2}) which contradicts Lemma \ref{2.26n}. Hence
also $\sigma_A$ satisfies (\ref{eq1.1}) and (\ref{eq1.2}). Let now
$\alpha > 0$ be arbitrary, choose mutually disjoint finite sets $A_k
\subseteq \bN$ so that for all $k \in \bN$ we have
\begin{equation*}
\alpha \le (\sum_{n \in A_k}\sigma_A(n)^r)^{\frac{1}{r}} \le 2\alpha
\end{equation*}
and put $\beta_k = (\sum_{n\in A_k}\rho_A(n)^r)^{\frac{1}{r}}$ for all
$k \in \bN$. Again Proposition \ref{thm1.1a} shows that
$X_p(\alpha,(\beta_k))$ is cb-isomorphic to a cb-complemented subspace
of $X_p(\sigma,\rho)$ and by choosing a subsequence of $(\beta_k)$
tending sufficiently fast to
$\beta = \liminf \beta_k >0$. we obtain that $X_p(\alpha, \beta)$
is cb-isomorphic to a cb-complemented subspace of $X_p(\sigma,
\rho)$. If $\beta =0$ we have $X_p(\alpha, \beta) = X_{p,c_p}(\alpha)$
and this violates (\ref{eq2.2ny}) of Proposition \ref{prop2.24n} for $\alpha$ small enough. If $\beta>0$,
then $\beta \le 2\delta \alpha$ and this violates (\ref{eq2.1n}) of
Proposition \ref{prop2.24n} for $\delta$ small enough. By choosing
$\alpha$ small enough (\ref{eq2.2n}) is violated and we have reached a
contradiction.

Interchanging the roles of $\sigma$ and $\rho$ in the
argument above we can conclude that there is a $\varepsilon >0$ and a
$\delta>0$ so that
\begin{eqnarray}
\sum_{\{\sigma_n \le \delta \rho_n, \sigma_n \le \varepsilon\}}\sigma_n^r < \infty
\\ \label{eq2.4n}
\sum_{\{\rho_n \le \delta \sigma_n, \rho_n \le \varepsilon\}}\rho_n^r <
\infty. \label{eq2.5n}
\end{eqnarray}

Let $A$ be as above and put

\begin{eqnarray*}
B & = & \{n \in \bN \mid \delta\rho_n < \sigma_n <
\frac{1}{\delta}\rho_n\}\\
D & = & \{n \in \bN \mid \sigma_n \le \delta \rho_n\}
\end{eqnarray*}
and define the sequences $(\sigma_A)$, $(\sigma_B)$, $(\sigma_D)$,
$(\rho_A)$, $(\rho_B)$, and $(\rho_D)$ as before. We can
then write
\begin{equation*}
X_p(\sigma, \rho) = X_p(\sigma_A, \rho_A) \oplus X_p(\sigma_B, \rho_B)
\oplus X_p(\sigma_D, \rho_D).
\end{equation*}

By Lemma \ref{2.25n} $X_p(\sigma_A, \rho_A) \oplus X_p(\sigma_D,
\rho_D)$ is cb-isomorphic to $l_p^N \oplus_p (C_p \cap R_p)^M\oplus_p C_p^K
\oplus_p R_p^L$ for some
$0\le k,L,M,N \le \infty$. $X_p(\sigma_B, \rho_B)$ is cb-somorphic to
$X_p(\sigma_B)$ and since $\sigma_B$ satisfies (\ref{eq1.1}) and
(\ref{eq1.2}) it contains cb-complemente copies of $l_p \oplus_p (C_p \cap
R_p)$ which shows that $X_p(\sigma,\rho)$ is cb-isomorphic to
$X_p(\sigma_B)\oplus_p C_p^K \oplus_p R_p^L$. This finishes the proof since
clearly $X_p(\sigma_B)$ is cb-isomorphic to $X_p(\sigma)$. Obviously the
$C_p$- and $R_p$-terms do not appear in case $\sigma_n \to 0$ and $\rho_n
\to 0$.
\eproof

%% file: rosmatrixnew11.tex
\section{Operator space properties of the matricial Rosenthal spaces}
\label{sec3}
\setcounter{equation}{0}
In this section we will discuss the operator space structure
of the matricial Rosenthal spaces. As before we let $p>2$, $\frac{1}{2}
= \frac{1}{p} + \frac{1}{r}$, and let
$\sigma$ be a sequence with $\sigma_n \ge 0$. $(\xi_n)$ denotes the unit
vector basis of $\ell_2$. Throughout the rest of the paper we let ${\mathcal R}$
denote the hyperfinite $II_1$ factor defined as the $\sigma$-weak
closure of the infinite tensor product $\otimes_{n\in \bN} M_2$ in the
GNS-construction with respect to the tracial trace
$\tau_{{\mathcal R}} = \otimes_{n \in \bN} \frac{tr}{2}$.

We start with the following result on $Y_p(\sigma)$:

\begin{prop} \label{cr}
$Y_p (\si)$ is complemented in $L_p({\mathcal R})$.
\end{prop}

{\bf Proof:} Let $\mu$ denote the Lebesgue measure on $]0,\infty[$ and
let $A_n\subset ]0,\8[$ be disjoint sets with $\mu(A_n) = \sigma_n^r$
for all $n \in \bN$. We consider the subspace
$V \subset L_p((0,\8);S_p) \cap L_2^{r_p\cap c_p}(]0,\8[;S_2)$
defined as the closure of $\{ \summ_n \mu(A_n)^{-\frac{1}{p}} 1_{A_n}x_n \pl |\pl
x_n \in S_p^n \}$.

Given $X_n\in S_p\ten S_p^n$, we have
\begin{equation*}
 \noo \summ_n \mu(A_n)^{-\frac{1}{p}} 1_{A_n} X_n\rrm_{L_p(S_p)}= \kla \summ_n \noo X_n\rrm_p^p \mer^{\frac1p}.
\end{equation*}
Further
\begin{eqnarray*}
  & &\noo \summ_n  \mu(A_n)^{-\frac{1}{p}} 1_{A_n} X_n\rrm_{S_p[L_2 ^{c_p}]}
  \lel    \noo \kla \summ_n \mu(A_n)^{1-\frac{2}{p}} (id\ten tr)(X_n^*X_n) \mer^{\frac12} \rrm_{S_p} \\
  & & \lel  \noo \kla \summ_n \si_n^2 (id\ten tr)(X_n^*X_n) \mer^{\frac12} \rrm_{S_p}
\end{eqnarray*}
The calculation   for the row term is similar. Comparing this with 
(\ref{eqyp}) we obtain that $V$ is
cb-isomorphic to $Y_p(\sigma)$.

For every $n\in \bN$ we let $p_n$ denote the orthogonal projection of
$\ell_2$ onto $\span\{\xi_n \mid \frac{n(n-1)}{2} + 1 \le k \le
\frac{n(n+1)}{2} \}$. Since $B=\{\sum_n 1_{A_n} \ten x_n \pl |
 \pl x_n=p_nx_np_n\}$ is a von Neumann
subalgebra of $L_\8((0,\8);B(\ell_2))$ and the restriction of the
trace is normal on $B$, we deduce from \cite{Ta} that there is
conditional expectation
\[ E(x) \lel \summ_n 1_{A_n} \ten \intt_{A_n} p_nx(t)p_n \frac{dt}{\mu(A_n)} \]
which is completely contractive on $L_p((0,\8);S_p)$ for all
$1\le p\le \8$. Clearly $E$ is  a projection onto $V$ and hence
$V$ is cb-complemented in $L_p((0,\8);S_p)\cap L_2^{c_p\cap
r_p}((0,\8);S_2)$.  According to \cite{J2} the latter space is
cb-isomorphic to $L_p({\mathcal R})$ and the assertion is
proved.\qed

\begin{rem} According to \cite{J2}, the spaces
\[ L_p((0,\8);S_p)\cap L_2^{c_p}((0,\8);S_2)\quad \mbox{and} \quad
  L_p((0,\8);S_p)\cap L_2^{r_p}((0,\8);S_2)\]
are cb-isomorphic to completely complemented subspaces in
$L_p({\mathcal R}\ten B(\ell_2))$ and hence the same argument as
above shows that $Y_{p,c_p}$ and $Y_{p,r_p}$ are cb-isomorphic to
cb-complemented subspaces of $L_p({\mathcal R}\ten B(\ell_2))$.
However, in general we cannot expect
a  cb-embedding into $L_p({\mathcal R})$. Indeed, from Theorem
\ref{thm1.3} it follows that if $\sigma$ satisfies (\ref{eq1.1}) and
(\ref{eq1.2}), then $S_p$ cb-embeds into $Z_p(\sigma)$ but it does not
embed into $L_p(\mathcal R)$ according to a result of Suckochev
\cite{Su}. Hence $Z_p(\sigma)$ does not cb-embed into $L_p(\mathcal R)$.
\end{rem}

\begin{cor}
\label{cor3.3}
 The spaces $Y_{p}(\si)$, $Y_{p,c_p}(\sigma)$ and
$Y_{p,r_p}(\sigma)$ have the
$\gamma_p$-AP.
\end{cor}

{\bf Proof:} Since $L_p({\mathcal R}\ten B(\ell_2))$ is the $L_p$
space of an injective von Neumann algebra, this space the
$\gamma_p$-AP. The  $\gamma_p$-AP passes to complemented
subspaces.\qed

We now turn our attention to the space $Z_p(\sigma)$ but for this we
need some preliminary results.

Let $m,n \in \bN$ and let $D$ be a positive $m\times m$ diagonal matrix
with $tr(D)=1$. We define $Z_p^m(n,D)$ to be the subspace of
$S_p^m\oplus_p C_P^{m^2} \oplus_p R_p^{m^2}$ defined by:
\begin{equation*}
Z_p^m(n,D) =
\{(x,n^{\frac{1}{r}}xD^{\frac{1}{r}},n^\frac{1}{r}D^\frac{1}{r}x) \mid
x\in S_P^m \}.
\end{equation*}
Here we consider $xD^\frac{1}{r}$ as an element of $C_p^m(C_p^m) =
C_p^{m^2}$ and $D^\frac{1}{r}x$ as an element of
$R_p^m(R_p^m)=R_p^{m^2}$. The spaces $Z_{p,c_p}^m(n,D)$ and $Z_{p,r_p}^m(n,D)$
are defined similarly as subspaces of $S_p^m \oplus_p C_p^{m^2}$,
respectively $S_p^m \oplus_p R_p^{m^2}$.

For every $1\le i \le n$ we define $\Psi_i:S_p^m \to S_p^{m^n} =
S_p^{\otimes_n}$ by
\begin{equation*}
\Psi_i(x) = D^\frac{1}{p} \otimes \cdots \otimes D^\frac{1}{p} \otimes
x \otimes D^\frac{1}{p} \otimes \cdots D^\frac{1}{p}
\end{equation*}
for all $x\in S_p^m$ where $x$ is the ith factor. Further we put
\begin{eqnarray*}
U_p(x) & = & n^{-\frac{1}{p}}\sum_{i=1}^n \varepsilon_i\Psi_i(x) \quad \mbox{for all $x\in S_p^m$}\\
U_{p,c}(x) & =& n^{-\frac{1}{p}}\sum_{i=1}^n \varepsilon_i\Psi_i(x)\otimes e_{i1} \quad
\mbox{for all $x\in S_p^m$}\\
U_{p,r}(x) & =& n^{-\frac{1}{p}}\sum_{i=1}^n \varepsilon_i\Psi_i(x)\otimes e_{1i} \quad
\mbox{for all $x\in S_p^m$}
\end{eqnarray*}
where $(\varepsilon_i)$ is the sequence of Rademacher functions on $[0,1]$.

\begin{theorem}
\label{thm3.1}
$U_p$ acts as a cb-isomorphism of $Z_p^m(n,D)$ onto its image which is
cb-complemented in $L_p([0,1],S_p^{m^n})$ with cb-norms only depending
on $p$.
Similarly $Z_{p,c_p}^m(n,D)$ ($Z_{p,r_p}^m(n,D)$) is cb-complemented in
$L_p([0,1];S_p^{m^n}\otimes C_p^n)$ ($L_p([0,1];S_p^{m^n}\otimes R_p^n)$)
via the map $U_{p,c}$ ($U_{p,r}$).
\end{theorem}
\bproof
Let $\{x_{jk} \mid 1\le j,k \le m \} \subseteq S_p$. If we for every
$1\le i \le n$ put
\begin{equation*}
Y_i = \varepsilon_i \sum_{j,k}x_jk \otimes \Psi_i(e_{jk}),
\end{equation*}
then the $Y_i$'s are independent in the sense of \cite{JXII} and have
mean zero.

Therefore, if we put $E(x\otimes y) = tr(D^{1-\frac{2}{p}}x)y$ for all
$x\in S_p^{m^n}$ and all $y\in S_p$ and let  ``$\sim$'' denote a two-sided inequality with constants
only depending on $p$, \cite[Theorem 1.2]{JXII} gives that
\begin{eqnarray}
 & & \|\sum_{j,k}x_{jk} \otimes
U_p(e_{jk})\|_{S_p[L_p(S_p^{m^n})]} =
n^{-\frac{1}{p}}\|\sum_{i=1}^n Y_i\|_{S_p[L_p(S_p^{m^n})] } \label{eq3.1} \\
& & \sim
n^{-\frac{1}{p}}\max\{(\sum_{i=1}^n\|Y_i\|_{S_p[L_p(S_p^{m^n})]}^p)^{\frac{1}{p}},
\|(\sum_{i=1}^nE(Y_i^*Y_i))^{\frac{1}{p}}\|_{S_p}, \|(\sum_{i=1}^n
E(Y_iY_i^*))^{\frac{1}{p}}\|_{S_p}\}.   \nonumber
\end{eqnarray}

For all $i \le n$ we easily get that
\begin{equation*}
\|Y_i\|_{S_p[L_p(S_p^{m^n})]} = \|\sum_{j,k} x_{jk}\otimes
e_{jk}\|_{S_p[S_p^m]}
\end{equation*}
Further
\begin{eqnarray*}
\|(\sum_{i=1}^n E(Y_i^*Y_i))^{\frac{1}{2}}\|_{S_p} & = &
n^{\frac{1}{2}}
\|(\sum_{j,k}\sigma_k^{1-\frac{2}{p}}x_{jk}^*x_{jk})^{\frac{1}{2}}\|_{S_p}
\lel  n^{\frac{1}{2}}\|\sum_{k=1}^m
\sigma_k^{\frac{1}{r}}\sum_{j=1}^m x_{jk} \otimes
e_{jk}\|_{S_p[C_p^{m^2}]}
\end{eqnarray*}

and similarly
\begin{equation*}
\|(\sum_{i=1}^n E(Y_iY_i^*))^{\frac{1}{2}}\|_{S_p} =
n^{\frac{1}{2}}\|\sum_{j=1}^n\sigma_j^{\frac{1}{r}}\sum_{k=1}^m x_{jk}
\otimes e_{jk}\|_{S_p[R_p^{m^2}]}.
\end{equation*}
Combining these calculations with (\ref{eq3.1}) we get that $U$ is a
cb-isomorphism of $Z_p^m(n,D)$ onto its image.

For every $1 \le i \le n$ we define $\Psi_i': S_{p'}^m \to S_{p'}^{m^n}$
by
\begin{equation*}
\Psi_i'(x) = D^{\frac{1}{p}}\otimes \cdots \otimes D^{\frac{1}{p}}
\otimes x \otimes D^{\frac{1}{p}} \otimes \cdots \otimes D^{\frac{1}{p}}
\end{equation*}
for every $x \in S_{p'}$ where $x$ is the ith factor and
$U_{p'} = \sum_{i=1}^n \varepsilon_i \Psi_i'(x)$ for all $x \in
S_{p'}$. Using \cite[Theorem 4.3]{JXII} we can in a similar manner as
above obtain that $U_{p'}$ acts as a cb-bounded operator from
$Z_p^m(n,D)^*$ to $L_{p'}([0,1],S_p^{m^n})$ It is readily verified that
$U_pU_{p'}^*$ is a cb-bounded projection of $L_p([0,1],S_p^{m^n})$ onto
the range of $U_p$.

The argument for $U_{p,c}$ and $U_{p,r}$ can be done similarly.
\eproof

We are now able to prove:

\begin{theorem}
\label{thmZ}
Let $2\le p,r< \8$ such that $\frac12=\frac1p+\frac1r$.
If $\si$ is a sequence of positive numbers such that $\si \notin
\ell_r$ and $\liminf_n \si_n=0$, then  the spaces $Y_p(\sigma)$,
$Y_{p,r_p}(\sigma)$, $Y_{p,c_p}(\sigma)$, $Z_p(\si)$,
$Z_{p,r}(\si)$, $Z_{p,c}(\si)$   are ${\mathcal COS}_p$ spaces.
\end{theorem}

{\bf Proof:}  Let us consider
 \[ s_j \lel \summ_{k=1}^j \si_k^r \]
By assumption $s_j$ tends to $\8$ and hence we can
find a subsequence $(j_k)$  and integers $n_k$ such
that
 \[ n_k \le s_{j_k} \le n_k+1 \pl .\]
By definition $Z_p$, $Z_{p,c}$, $Z_{p,r}$  is  the
closure of $\bigcup_k Z_p^{j_k}$,  $\bigcup_k
Z_{p,c}^{j_k}$, $\bigcup_k Z_{p,r}^{j_k}$,
respectively. Fix $k\in \nz$ and define $\rho_k \lel
s_{j_k}^{-1}(\si_j^r)_{j\le j_k}$. The map

 \[ w(x) \lel (x, n_k^{\frac{1}{r}}xD_{\rho_k}^{\frac1r}, n_k^{\frac1r}D_{\rho_k}^{\frac1r}x) \pl  \]
yields an isomorphism between  $Z_p^{j_k}(\si)$ and $Z_p(n_k,D_{\rho_k})$.
Indeed, for $\si_k=(\si_{j})_{j\le j_k}$ we have
 \[ n_k^{\frac1r} D_{\rho_k}^{\frac1r} \lel  \kla \frac{n_k}{s_{j_k}}\mer^{\frac1r}  D_{\si_k} \]
and
 \[ 1 \kl  \kla \frac{n_k}{s_{j_k}}\mer^{\frac1r} \kl (1+\frac{1}{n_k})^{\frac1r} \kl 2\pl. \]
Hence by Theorem \ref{thm3.1} $Z_p^{j_k}(\si)$ has the $\gamma_p$--AP
with a constant only depending $\sigma$ and $p$ and therefore $Z_p(\sigma)$ has
the $\gamma_p$--AP. Similarly for $Z_{p,c_p}(\si)$ and
$Z_{p,r_p}(\si)$. $Y_p(\sigma)$, $Y_{p,c_p}(\sigma)$, and
$Y_{p,r_p}(\sigma)$ have the $\gamma_p$--AP by Corollary
\ref{cor3.3}. Since $\liminf_n \si_n=0$,
we can find a subsequence $\si'=\si_{n_k}$ such
that $(\si_{n_k})\in \ell_r$. Then the map
$M_r:S_p\to C_p(\nz^2)$  defined by
$M_r(x)=xD_{\si'}$ is completely bounded and
similarly, $M_l:S_p\to R_p(\nz^2)$ defined by
$L_l(x)=D_{\si'}x$ is completely bounded. If
$A=\{n_k:k\in \nz\}$, then the subspace
$Z_A=\{(x_{ij})\p |\p i\in A, j\in A\}$ is
cb-isomorphic to $S_p$ and is complemented in
$Z_{p}(\si)$, $Z_{p,c}(\si)$, and $Z_{p,r}(\si)$,
respectively. By the definition of $Y_{p}(\si)$ we
deduce that $Y_A=\{(x_k)_{k} \p |\p k\in A, x_k\in
M_{m_k}\}$ is cb-isomorphic to $(\sum_{k\in A}
\oplus_p S_p^{m_k})_p$ and cb-complemented. Thus all
these spaces contain $S_p^n$'s uniformly
complemented. According to \cite[Theorem 2.2]{JNRX},
we deduce that these spaces are ${\mathcal COS}_p$
spaces.\qed

%% file: uncomp2.tex
\section{Uncomplemented copies of some $\cO\cL_p$--spaces}
\label{sec4}
\setcounter{equation}{0}

Throughout this section, $2 < p < \infty$, unless specified
otherwise.

\begin{theorem}
\label{thm4.1}

Let $X$ and $Y$ be subspaces of rectangular
$\cO\cL_p$ spaces so that $X$ is completely isomorphic
to a subspace of $Y$. Then $\ell_p(Y)$ (respectively, $S_p[Y]$)
contains an uncomplemented completely isomorphic copy of
$\ell_p(X)$ (respectively, $S_p[X]$).
\end{theorem}

Before proving the theorem, we formulate a corollary of it.

\begin{cor}
\label{cor4.2}
\begin{itemize}
\item[(a)] Suppose $X$ is one of the following operator spaces: $\ell_p$,
$S_p$, ${\K}_p$, or $L_p({\mathcal R})$. Then $X$ contains an uncomplemented copy of itself.
\item[(b)] Suppose $\N$ is a group von Neumann algebra with QWEP,
and $X$ is either $\ell_p(L_p(\N))$, or $S_p[L_p(\N)]$.
Then $X$ contains an uncomplemented copy of itself.
\end{itemize}
\end{cor}

\bproof
All the spaces listed in parts (a) and (b) are $\cO \cL_p$-- spaces
(see \cite{JR} for the spaces from part (b)).
Moreover, any of the spaces $X$ listed in part (a) is completely
isomorphic to $\ell_p(X)$, by Pe{\l}czy\'{n}ski's decomposition method.
The same argument shows that for $\N$ as in part (b) $S_p[L_p(\N)]$
is completely isomorphic to $\ell_p(S_p[L_p(\N)])$.
\eproof

To establish Theorem \ref{thm4.1}, consider a finite dimensional
version of the Rosenthal space. More precisely, if
$\sigma = (\sigma_n)_{n \in \bN}$ is a sequence of positive numbers,
then we let $X_p^m(\sigma)$ be the linear span of the first $m$ vectors of the canonical basis of $X_p(\sigma)$. By Corollary \ref{njn2.2} 
there exists $\lambda > 0$, and a sequence
$(k_m)_{m \in \N}$, s.t. $\ell_p^{k_m}$ contains a
$\lambda$-completely complemented $\lambda$--completely isomorphic
copy of $X_p^m(\sigma)$.

Now suppose the sequence $(\sigma_n)$ satisfies (1.5) and (1.6).
By \cite{R}, if $P_m$ is a projection from
$\ell_p^m \oplus_p R_p^m \oplus_p C_p^m$
onto the ``natural'' copy of $X_p^m(\sigma)$, then
$\lim_m \|P_m\| = \infty$. By \cite{Lu} (see also \cite{P1}),
$\ell_p^m \oplus_p R_p^m \oplus_p C_p^m$ embeds into
$\ell_p^{3^m}$ $c_p$--completely isomorphically.
Thus, there exists a sequence $(T_m)$ of complete contractions
$T_m : X_p^m(\sigma) \to \ell_p^{3^m}$ so that
$\|T_m^{-1}\|_{cb} \leq c_p$, and $\lim_m \|Q_m\| = \infty$
whenever $Q_m$ is a projection from $\ell_p^{3^m}$ onto $range (T_m)$.

The properties of the spaces $X_p^m(\sigma)$ yield:

\begin{lemma}
\label{lemma4.3}
$\ell_p$ contains an uncomplemented completely isomorphic copy
of itself.
\end{lemma}
\bproof
Suppose the sequence $(\sigma_m)$ satisfies (1.5) and (1.6).
Consider the spaces $Y = (\sum_m \ell_p^{3^m})_p$, and
$Z = (\sum_m T_m(X_p^m(\sigma)))_p$. By the discussion preceding
the statement of this lemma, $Z$ is an uncomplemented subspace
of $Y$. Moreover, $Y$ is completely isometric to $\ell_p$.
It remains to show that $Z$ is completely isomorphic to $\ell_p$.
To this end, note that $Z$ is completely isomorphic to a
completely complemented subspace of
$(\sum_m \ell_p^{k_m})_p \sim \ell_p$.
Moreover, $Y$ contains a completely complemented copy of $\ell_p$.
As $\ell_p = \ell_p(\ell_p)$, we complete the proof by applying
a Pe{\l}czy\'{n}ski decomposition method.
\eproof

We need yet another lemma.

\begin{lemma}
\label{lemma4.5}
Suppose $X$ is a rectangular $\cO\cL_p$ space, and $T$
is a complete isomorphism from $\ell_p$ onto a subspace.
Then $T \otimes I_X$ is a complete isomorphism from $\ell_p(X)$
onto its range, viewed as a subspace of $\ell_p(X)$.
\end{lemma}

\bproof
We can assume that $T$ is a complete contraction and let
$c = \|T^{-1}\|_{cb}$. It suffices to show that
$T \otimes I_{S_p^N} : \ell_p(S_p^N) \to \ell_p(S_p^N)$ is a complete
contraction, and $\|(T \otimes I_{S_p^N})^{-1}\|_{cb} \leq c$.
To complete the proof identify $\ell_p(S_p^N)$ with
$S_p^N[\ell_p]$ and apply Proposition \ref{prop01}.
\eproof

\begin{rem}
The same result also holds for complete isomorphisms from $S_p$
onto its subspaces.
\end{rem}

\noindent {\bf Proof of Theorem \ref{thm4.1}:}
Suppose $X$ and $Y$ are subspaces of rectangular $\cO\cL_p$--spaces
and $S : X \to Y$ is a complete isomorphism. Let $T : \ell_p \to \ell_p$
be a complete isomorphism with an uncomplemented range (such a $T$
exists, by Lemma \ref{lemma4.3}). By Lemma \ref{lemma4.5} $T \otimes S$
determines a complete isomorphism from $\ell_p(X)$ onto a
subspace of $\ell_p(Y)$. It remains to show that
$range(T \otimes S)$ is uncomplemented. Indeed, suppose
for the sake of contradiction that there exists a projection
$P$ from $\ell_p(Y)$ onto $range (T \otimes S)$.
Pick $x \in X \backslash \{0\}$ and denote by $Q$ a bounded
projection onto $\span (Sx)$. As $T$ is a complete isomorphism,
$\tilde{Q} = I_{range (T} \otimes Q$ is a completely bounded
projection from $range (T \otimes S)$ onto $range (T) \otimes \span (Sx)$.
Hence $\tilde{Q} \circ P|_{\ell_p \otimes \span (Sx)}$ is a bounded projection
from $\ell_p \otimes \span (Sx)$ onto $range (T) \otimes \span (Sx)$ which contradicts
the fact that $range (T)$ is uncomplemented.
\eproof
\begin{cor}
\label{cor4.6}
Suppose $\N$ is a von Neumann algebra equipped with a normal
semi-finite faithful trace which is not of type $I$. Then
there exists an uncomplemented subspace $X$ of $L_p(\N)$ completely
isomorphic to $L_p({\mathcal R})$
\end{cor}
\bproof
By \cite{JNRX} (see also \cite{Ma}) $L_p(\N)$ contains a 
(completely contractively complemented) subspace $Y$,
completely isometric to $L_p({\mathcal R})$. By Theorem \ref{thm4.1} $Y$ contains an
uncomplemented copy of $L_p({\mathcal R})$.
\eproof

\begin{cor}
\label{cor4.7}
\begin{itemize}
\item[(1)] Every infinite dimensional rectangular $\cO\cL_p$--space contains
an uncomplemented copy of $\ell_p$.
\item[(2)] Every infinite dimensional $\cO\cS_p$--space contains an
uncomplemented copy of $(\sum_n S_p^n)_p$.
\end{itemize}
\end{cor}
\bproof
By \cite{JNRX} any $\cO \cL_p$--space $X$ (with $1 < p < \infty$) embeds
completely isometrically (and even completely contractively
complementedly) into $\Pi_{\cU} S_p$, where $\cU$ is an
ultrafilter. By \cite{Re} and \cite{RX}
$X$ contains a completely isomorphic (and even completely
complemented) subspace $Y$, completely isomorphic to $\ell_p$.
Moreover, if $X$ is an $\cO\cS_p$--space, then it contains a
subspace $Y$, completely isomorphic to $(\sum_n S_p^n)_p$.
In either case an application of Theorem \ref{thm4.1} completes the proof.
\eproof

%% file: refnew11.tex
\vspace{1cm}

\noindent Department of Mathematics,\\ University of Illinois at 
Urbana--Champaign,\\ 1409 W. Green Street,\\ Urbana, Il 61801, U.S.A.\\
junge@.math.uiuc.edu\\

\noindent Department of Mathematics and Computer Science,\\ University of 
 Southern Denmark,\\ Campusvej 55, DK-5230 Odense M, Denmark,\\
njn@imada.sdu.dk\\

\noindent Department of Mathematics, \\ University of California, Irvine, 
\\ 103 MSTB, Irvine, CA 92697--3875, U.S.A.\\
toikhber@math.uci.edu